\documentclass[preprint,12pt]{elsarticle}
\usepackage[margin=2cm]{geometry}
\usepackage{lineno,hyperref}
\modulolinenumbers[5]
\usepackage{xcolor}
\usepackage{amssymb,amsmath,color,graphicx, amsbsy}
\allowdisplaybreaks
\journal{International Journal of Electrical Power and Energy Systems}

\begin{document}

\begin{frontmatter}



\title{Strategic Bidding in Electricity Markets with Convexified AC Market-Clearing Process}


\author[inst1]{Arash Farokhi Soofi}

\affiliation[inst1]{San Diego State University, United States},%Department and Organization
            % addressline={Address One}, 
            % city={City One},
            % postcode={00000}, 
            % state={State One},
            % country={Country One}}

\author[inst1]{Saeed D. Manshadi\corref{mycorrespondingauthor}}
% \author[inst1,inst2]{Author Three}

% \affiliation[inst2]{organization={Department Two},%Department and Organization
%             addressline={Address Two}, 
%             city={City Two},
%             postcode={22222}, 
%             state={State Two},
%             country={Country Two}}

\begin{abstract}
This paper presents a framework to solve the strategic bidding problem of participants in an electricity market cleared by employing the full AC Optimal Power Flow (ACOPF) problem formulation. Traditionally, the independent system operators (ISOs) leveraged DC Optimal Power Flow (DCOPF) problem formulation to settle the electricity market. The main quest of this work is to find \textit{ what would be the challenges and opportunities if ISOs leverage the full ACOPF as the market-clearing Problem (MCP)?% instead of the current practice of using the DC form?
} 
This paper presents tractable mathematical programming with equilibrium constraints for the convexified AC market-clearing problem. Market participants maximize their profit via strategic bidding while considering the reactive power dispatch of generation units. The equilibrium constraints are procured by presenting the dual form of the relaxed ACOPF problem. The strategic bidding problem with ACOPF-based MCP improves the exactness of the location marginal prices (LMPs) and profit of market participants compared to the one with DCOPF. It is shown that the strategic bidding problem with DCOFP-based MCP is unable to model the limitation of reactive power support. The presented results display cases where the proposed strategic bidding method renders $52.3\%$ more  profit for the Generation Company (GENCO) than the DCOPF-based MCP model. The proposed strategic bidding framework also addresses the challenges in coupling real and reactive power dispatch of generation constraints, ramping constraints, demand response implications with curtailable and time shiftable loads, and AC line flow constraints. Therefore, the presented method will help market participants leverage the more accurate ACOPF model in the strategic bidding problem. 
\end{abstract}

\begin{graphicalabstract}
\includegraphics[width=\linewidth]{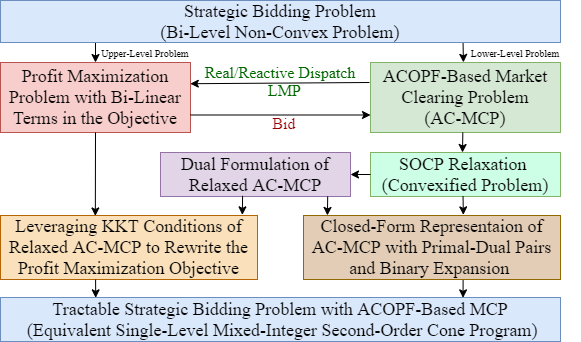}
\end{graphicalabstract}

\begin{highlights}
% \begin{itemize}
 \item Utilize ACOPF-based market-clearing problem with real and reactive power dispatches.
    \item Present a tractable closed-form for the ACOPF-based market-clearing problem.
    \item Show the limitations of the DCOPF-based market-clearing compared to the AC form.
    \item Investigate various market implications of employing ACOPF-based market-clearing.
% \end{itemize}
\end{highlights}

\begin{keyword}
 electricity market \sep strategic bidding \sep convex relaxation \sep second-order cone programming \sep bi-level optimization
\end{keyword}

\end{frontmatter}


\section{NOMENCLATURE}
\label{sec:Nom}

\subsection*{Variables}
\noindent \begin{tabular}{ l p{16cm} } \allowdisplaybreaks

$C^{'u,t}_g$ & Bid of segment $u$ of unit $g$ at time $t$  \\ %in $\$/MWh$

$C_{l,P}^{'z,t}$ & Bid of segment $z$ of load $l$ at time $t$\\

$c,s$  & Lifting operator terms for SOCP relaxation\\

$e_i^t$, $f_i^t$ & Real and imaginary parts of voltage phasor of bus $i$ at time $t$\\

$I_g^{c,u,t}$ & Binary variable representing the bidding of segment $u$ of generation unit $g$ at time $t$\\

% $\gamma_g^{s}$ & Bidding strategy of segment $s$ of generation unit $g$ \\

% $V_{i}^t$ & Voltage phasor of bus $i$ at time $t$\\

% $f_i$ & Imaginary part of voltage phasor of bus $i$ \\

$P_g^t$, $Q_g^t$ & Real and reactive power dispatch of generation unit $g$ at time $t$\\

$P_d^{z,t}$ & Dispatch to segment $z$ of load $d$ at time $t$\\

$P_g^{u,t}$ &  Real power dispatch of segment $u$ of generation unit $g$ at time $t$ \\ 
$P_l^{s,t}$, $Q_l^{s,t}$ &  Sending real and reactive power of line $l$ at time $t$\\

$P_l^{r,t}$, $Q_l^{r,t}$ &  Receiving real/reactive power from line $l$ at time $t$\\
%\end{tabular}
%\newpage
%\noindent \begin{tabular}{ ll p{6.7cm} } 

% $S_g^{u,t}$ & apparent power dispatch of segment $u$ of generation unit $g$\\
% $\theta_{i}^t$ &  Voltage angle of bus $i$ at time $t$\\

% $\theta_{ij}^t$ &  Phase angle difference of bus $i$ and bus $j$ at time $t$\\

$\mathcal{\mu}$ & Dual variable for inequality constraints\\

$\mathcal{\lambda}$ & %\hspace{0.7cm} 
Dual variable for equality constraints

% $\mathcal{\mu}$ & Dual variable corresponding to inequality constraints

\end{tabular}
\subsection*{Parameters}
\noindent \begin{tabular}{ l p{16cm} }\allowdisplaybreaks
$b_{l}^{sh}$  & The shunt susceptance of line $l$ \\

$g_{l},b_{l}$  & The real and imaginary parts of admittance of line $l$\\
%$G_{ii}$  & Shunt conductance at bus i\\
%$B$  & Susceptance matrix\\
%$b_{l}$  &  part of admittance of line $l$\\

$C_g^{P,u}$ &  Generation cost of segment $u$ of generation unit $g$ \\

$C_g^Q$ &  Reactive power generation cost of generation unit $g$ \\

$\underline{P}_g,\overline{P}_g$  & Minimum and maximum capacity of real power of generation unit $g$. \\

$\underline{Q}_g,\overline{Q}_g$  & Minimum and maximum capacity of reactive power of generation unit $g$\\

$\overline{P}_g^u$  & Maximum generation capacity of segment $u$ of generation unit $g$\\

$R_g^u, R_g^d$ & Upper/lower bounds of ramping limits of unit $g$\\

$p_d^t$, $q_d^t$ & Real and reactive demand $d$ at time $t$\\

$\overline{S_l}$  & Maximum apparent power capacity of line $l$\\
% $\overline{S}_g^u$  & Maximum apparent power of segment $u$ of generator $g$\\

$W_d^{z,t}$ & Willingness to pay of segment $z$ of demand $d$ at time $t$\\

$\overline{V}_i,\underline{V}_i$  & Maximum and minimum voltage magnitude at bus $i$\\

$\alpha_c$ & The coefficient of bidding strategy \\ 
$\underline{\alpha},\overline{\alpha}$ & The lower and upper limits of submitted bids\\
% \end{tabular}
% \newpage
% \noindent \begin{tabular}{ l p{7.5cm} }
% $\varrho_g$ & The profit of generation unit $g$\\
$\rho_g^{-/+,u}$ & Negative/positive slope of segment $u$ of piece-wise linear model of the coupling of $P_g$ and $Q_g$\\
% capability curve of generation unit $g$

% $\rho_l^{-/+,y}$ & Slope of segment $y$ of piece-wise linear function representing the coupling of $P_l$ and $Q_l$\\

$\chi_d^z$ & The portion of segment $z$ from load $d$
%$\overline{S_{ij}}$ & Maximum capacity of branch $ij$\\

%$\underline{S_{ij}}$ & Minimum capacity of branch ij\\

%$c^{B}[y]$ & Installation cost a battery unit at year $y$  \\

%$c^{R}[y]$ & Installation cost of an inverter  at year $y$  \\

%$I_{b}^{B}[y]$ & A binary variable indicating if battery storage unit $b$ is first installed at year $y$\\

%$I_{br}^{R}[y]$ & A binary variable indicating if inverter $r$ is first installed at year $y$\\

\end{tabular}
%\vspace{-1cm}
\subsection*{Sets}

\noindent \begin{tabular}{ ll p{12cm} }\allowdisplaybreaks

% $\mathcal{BP}$ & Set of all buspairs \\

% $\mathcal{BP}_{ij}$ & Set of the buspair between bus $i$ and bus $j$ \\

% $\mathcal{BP}_l$ & Set of the buspair connected by line $l$ \\

$\mathcal{BPF}_i$ & \hspace{0.15cm}Set of all buspairs originated from bus i\\

$\mathcal{BPT}_i$ & \hspace{0.15cm}Set of all buspairs destined to bus i\\

$\mathcal{D}, \mathcal{D}_c$ & \hspace{0.15cm}Sets of all loads and curtailable loads\\

$\mathcal{D}_i$ & \hspace{0.15cm}Sets of all loads connected to bus $i$\\

$\mathcal{F}_{bp} $ & \hspace{0.15cm}Set of buses which buspair $bp$ originated from \\

$\mathcal{F}_l $ & \hspace{0.15cm}Set of buses which line $l$ is originated from \\

$\mathcal{G}$ & \hspace{0.15cm}Set of all generation units  \\

$\mathcal{G}_B$ & \hspace{0.15cm}Set of all generation units belong to GENCO $B$ \\

$\mathcal{G}_{-B}$ & \hspace{0.15cm}Set of all generation units don't belong to GENCO $B$ \\

$\mathcal{G}_i$ & \hspace{0.15cm}Set of all generation units connected to bus $i$ \\

$\mathcal{I}_g$ & \hspace{0.15cm}Set of buses connected to  generation unit $g$ \\

$\mathcal{I}_d $ & \hspace{0.15cm}Set of buses connected to load $d$\\

$\mathcal{L}$ & \hspace{0.15cm}Set of all lines  \\

% $\mathcal{L}_{bp} $ & Set of all lines connecting buspair $bp$\\

$\mathcal{LF}_i$ & \hspace{0.15cm}Set of all lines originated from bus i\\

$\mathcal{LT}_i$ & \hspace{0.15cm}Set of all lines destined to bus i\\

$\mathcal{T}$ & \hspace{0.15cm}Set of time horizon\\
 
$\mathcal{T}_l $ & \hspace{0.15cm}Set of buses which line $l$ destined to \\

$\mathcal{T}_{bp} $ & \hspace{0.15cm}Set of buses which buspair $bp$ destined to \\

$\mathcal{N}$ & \hspace{0.15cm}Set of all buses \\

$\mathcal{U}_g$ & \hspace{0.15cm}Set of segments of generation unit $g$ \\
$\mathcal{Z}_d$ & \hspace{0.15cm}Set of segments of load $d$ 
\end{tabular}

\section{Introduction}
\label{sec:sample:appendix}
Before deregulation, the electricity industry was regulated by the federal energy regulatory commission and state public utility commissions. After the deregulation, electricity is evolved into a distributed commodity. The electricity market derives its price aiming to reduce the total cost of the network through the increase in competitiveness \citep{gao2015optimal}.
 In this paradigm, each market participant (e.g., GENCOs) wants to maximize its profit by leveraging an optimal bidding strategy in a bilateral electricity market consisting of GENCOs and loads \cite{song2002nash} or in a transmission-constrained network where GENCOs have incomplete information \cite{li2005strategic}. In electrical power systems, market participants submit their bids to the ISO. Then, the ISO will clear the market to determine the Locational Marginal Pricing (LMP) and the awarded hourly generation dispatch of participants over $24$ hours.
 The strategic bidding of GENCOs in electricity markets with the DC market-clearing Problem (DC-MCP) is extensively studied in the \color{black} literature. Each GENCO maximizes its profit given the DCOPF problem solved by ISO to clear the market. For instance, in \cite{ruiz2009pool}, the authors proposed procedure to derive strategic offers relies on a bi-level programming model whose upper-level problem represents the profit maximization of the strategic GENCOs. In contrast, the lower-level one represents the MCP and the corresponding price formation. The bi-level model is reduced to a mixed-integer linear programming problem using the strong duality theorem and the Karush–Kuhn–Tucker (KKT) optimality conditions. In \cite{pozo2011finding}, authors proposed a compact formulation to find all pure Nash equilibria with stochastic demands based on the Stackelberg game. In \cite{baslis2011mid}, the authors proposed a yearly stochastic self-scheduling model for a price-maker hydro producer. In \cite{kazempour2011strategic}, a model aimed at helping strategic producers in making informed decisions on generation capacity investment is proposed. 
In all these articles, the ISO minimizes the total cost of operation subjected to physical constraints of the power system. 
Thus, the strategic bidding problem is formulated as a \textit{bi-level} optimization problem that can be reformulated as a single \textit{Mathematical Program with Equilibrium Constraints} (MPEC) with DC approximation power ﬂow constraints \cite{hobbs2000strategic}. 
  \par
  With the recent advancement to find a polynomial-time solution for the full ACOPF problem, ISOs will soon consider clearing the market based on ACOPF problem formulation. This paper envisions a strategic bidding problem based on the assumption that ISOs will adopt ACOPF problem formulation for their market-clearing process. 
 The AC power flow constraints and the network losses are not directly considered in the market-clearing based on the DCOPF problem. Thus, the MCP might not present a realistic representation of the system. There are two major approaches to model network losses in the DCOPF calculation. In marginal loss modeling, the loss of lines is presented as a linear function of nodal injections.
 In \cite{eldridge2017marginal}, the authors investigate some aspects of including a marginal line loss approximation in the DCOPF. The impact of the loss of lines on the LMP of buses is illustrated in \cite{fu2006different}. A matrix loss distribution framework is developed in \cite{sarkar2009dcopf} to gain more design flexibility in the distribution of system losses. Each line is considered separately for loss distribution. Another approach is to model the loss of lines in the DCOPF problem as a piece-wise linear function of the line flow, as shown in ~\cite{dos2010dynamic}. In \cite{akinbode2013fictitious}, authors examine the effect of non-physical losses on optimal power flow solutions and LMPs. Non-physical losses are created when segments of piece-wise linear functions used to approximate real power losses for a DCOPF problem are selected in the wrong order. A better curve fitting technique is presented in \cite{vaishya2019accurate} that can improve the power flow accuracy of the piece-wise linear loss modeling in the DCOPF calculation. These methods approximate the loss of the network with a linear or piece-wise linear function. Therefore the solution procured by the lower-level problem may be different from the solution of the ACOPF problem.

Different relaxation methods including SOCP \cite{jabr2008optimal}, Semi-Definite Programming (SDP) \cite{lavaei2011zero} and moment relaxation methods\cite{molzahn2014moment}, Quadratic Convex (QC) \cite{coffrin2015qc} presented a convex relaxed form of the full ACOPF problem with improved relaxation gap. Besides, the SOCP relaxation is strengthened with polyhedral envelopes \cite{kocuk2016strong}, McCormick envelopes \cite{bynum2018strengthened}, mixed SOCP/SDP moment relaxations \cite{molzahn2015mixed}, and cycle constraints \cite{soofi2020socp}.
 These developments motivated this research to utilize the ACOPF-based MCP for the strategic bidding problem by leveraging the dual form of the convexified ACOPF problem and procuring the MPEC problem with convexified AC constraints. Since the solution procured by the DCOPF problem formulation for the MCP may encounter inaccuracies, in this paper, the full ACOPF problem formulation is employed to increase the accuracy. Besides, the reactive power dispatch and flow in the electricity network cannot be considered in DCOPF-based MCP formulation.
Thus, the literature did not consider the limitation of reactive power support of generation units, the impact of the reactive power flow on the voltage magnitude of buses, and the capacity of lines. This paper aims to fill this gap by introducing a new formulation for the strategic bidding problem that maximizes the profit of real and reactive power dispatch of market participants (e.g., GENCOs) subjected to the ACOPF-based MCP.
\par
Leveraging the ACOPF problem formulation as the MCP of strategic bidding problem, instead of the DCOPF problem formulation, brings the following questions: \textit{Can we formulate strategic bidding problem once ISOs adopt ACOPF for their market-clearing prices, and if so, what are the challenges and opportunities? What would be market inefficiencies in rendering the LMPs under DC constraints once compared with a market-clearing process with AC constraints? What would be the impact of the limitation of reactive power support of generation units, demand-side management, coupling of real and reactive power dispatch of generation units, and ramping limits of generation units on the solution of the strategic bidding problem? What would be the impact of considering the reactive-LMP on the strategic bidding problem?}

\begin{figure}[h!]
\centering%width=\columnwidth
{\includegraphics[width=0.6\linewidth]{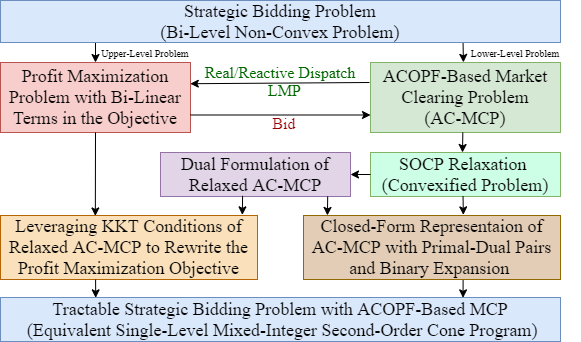}}
\caption{Transforming the non-convex strategic bidding problem with ACOPF-based MCP into a tractable problem }
\label{fig:flowchart}
\end{figure}

The main contributions of this paper are listed as follows.
\begin{itemize}

\item Propose a strategic bidding problem with ACOPF-based MCP for market participants to simultaneously consider real and reactive power dispatch. 

\item Procure the closed-form representation of the relaxed ACOPF-based MCP leveraging its dual formulation considering voltage limits as well as real and reactive power flow constraints.

\item Illustrate the capability of ACOPF-based MCP to capture the reactive power dispatch of generation units as well as voltage limitations and demonstrate the limitations of DCOPF-based MCP to enforce those constraints. 

\item Investigate the market implications of employing ACOPF problem formulation for the MCP, including reactive power support of generation units, demand-side management, coupling of real and reactive power dispatch of generation units, and ramping limits of generation units.
\end{itemize}
 The investigation procedure of this paper is presented in Fig. \ref{fig:flowchart}.
 The main goal of this procedure is to solve the strategic bidding problem with ACOPF-based MCP. As presented in Section 2, this problem is non-convex due to the non-convexity of ACOPF-based MCP, and non-linear due to the bi-linear terms multiplying LMPs with awarded dispatches in the objective. Therefore, the strategic bidding problem is formulated as a bi-level non-convex problem. The lower-level problem is the ACOPF-based MCP which passes the LMPs and awarded real and reactive power dispatches to the profit maximization problem. The SOCP relaxation method is deployed to relax the non-convexity of the ACOPF-based MCP in section 3.2. Once the dual form of the relaxed problem is procured in Section 3.3, the MCP is represented in a tractable mixed-integer conic closed-form with a set of primal, dual, the equality of primal and dual objectives constraints, and binary expansion in section 3.4. The upper-level problem aims to maximize the profit for market participants. Its bi-linear terms are replaced by the equivalent terms procured by leveraging the KKT condition of the relaxed ACOPF-based MCP as well as constraints in the dual form of AC-MCP as presented in Section 3.5. Thus, the non-convex bi-level  problem is transformed into a tractable single-level problem that can be solved with the off-the-shelf mixed-integer conic solvers.   

\section{Problem Formulation} \label{prob}
%\subsection{}
The strategic bidding problem for a GENCO, as one type of market participant, is formulated in \eqref{opf}, when ACOPF problem formulation is employed as the MCP. The presented problem is a bi-level optimization problem, where the upper-level problem is the profit maximization problem. The lower-level problem is the MCP. The objective of the upper-level problem is to maximize the profit of a GENCO, as shown in \eqref{opfUL}. The revenue is procured both by the real and reactive power generation of the unit that belongs to the GENCO, where the power is multiplied by the associated LMP of the connected bus. Besides, the GENCO expenses are presented by the summation of dispatches multiplied by their generation cost. More reactive power dispatch of the generation units leads to more real power loss. Besides, reactive power causes undesirable mechanical stresses and vibrations to the torque of the generation unit and may cause excessive wear and tear, which will decrease the generation unit's lifespan. Thus, the cost of dispatching reactive power of generation units is considered here. The upper bound and lower bound of the bid submitted by generation unit $g$ at time $t$ is presented in \eqref{opf_lim_bid}. The MCP with full ACOPF formulation in rectangular form is presented in \eqref{opfLL}-\eqref{opflimsend}. Here, the indicated variables after each colon in \eqref{opfpseg}-\eqref{opflimsend} introduce the dual variables corresponding to each constraint in the primal form of the MCP.
In the objective function of the MCP shown in \eqref{opfLL}, the total cost of bidding and non-bidding generation units including cost of real and reactive power dispatch is minimized. At the same time, the served demand is maximized based on the willingness to pay (WTP) of consumers. The total real power dispatch of generation unit $g$ at time $t$ is equal to the summation of the dispatch of each segment of generation unit $g$ at time $t$, as shown in \eqref{opfpseg}. The ramping up and down constraints of power generation unit $g$ are presented in \eqref{opf_ramp_t} and \eqref{opf_ramp_1}, respectively. The physical limits of real power generation of unit $g$ and real power generation of segment $u$ of unit $g$ at time $t$ are presented in \eqref{opfpg} and \eqref{opflimpseg}, respectively. The capability curve of synchronous generation units defines the boundaries within which it can deliver reactive power continuously without overheating \citep{nilsson1994synchronous}. The approximated piece-wise linear constraint representing the capability curve of generation units is presented in \eqref{opfqg_p}, where the upper bound and lower bound of the reactive power of unit $g$ at time $t$ is a function of the maximum/minimum reactive power dispatch of unit $g$, the slope of piece-wise linearized capability curve, and the real power generation dispatch of segment $u$ of unit $g$ at time $t$. It is assumed here that each generation unit only works in over-excitation or under-excitation modes over period $\mathcal{T}$. Since $\rho_g^{-,u}$ is negative, the upper bound of the reactive power of unit $g$ at time $t$ decreases when the real power of the generation unit increases. Since $\rho_g^{+,u}$ is positive, the lower bound of the reactive power of unit $g$ at time $t$ increases when the real power of the generation unit increases. The upper bound of the delivered real power to segment $z$ of fixed load $d$ at time $t$ is presented in \eqref{willingness_p}. Note that the summation of all $\chi_d^z$ over all segments of the fixed load is one. Thus, the delivered real power to fixed load $d$ at time $t$ is less than or equal to the demand of load $d$ at time $t$. 

The equality constraint given in \eqref{opf_fixed} represents the model of fixed load $d$ at time $t$ that should be served. The constraints modeling the shiftable and non-shiftable portions of curtailable loads are presented in \eqref{opf_flex}. The summation of served demand of segment $z$ of the shiftable curtailable load $d$ over $24$ hours is equal to the $\chi_d^z$ multiplied by the summation of demand of shiftable curtailable load $d$ over $24$ hours. This enables the ISO to serve the shiftable portion of the load during off-peak hours. Here, the curtailable load refers to a load that a portion of its demand can be curtailed at each hour, while the shiftable load is a load that the summation of its served demand over $24$ hours is equal to the summation of its demand over $24$ hours. %The served demand of segment $z$ of non-shiftable curtailable load $d$ at time $t$ is equal to the $\chi_d^z$ multiplied by the demand of non-shiftable curtailable load $d$ at time $t$.
The nodal balance equations for the real and reactive power at each bus are shown in \eqref{opfpbus} and \eqref{opfqbus}, respectively. The dual variables corresponding to constraints \eqref{opfpbus} and \eqref{opfqbus} are the LMPs for the real and reactive power at bus $i$. % and their unit are $\$/MWh$, $\$/MVARh$ respectively.
The real and reactive power flow sending through line $l$ are given in \eqref{opfpsend} and  \eqref{opfqsend}, respectively. The real and reactive power flow receiving from line $l$ are given in \eqref{opfprec} and \eqref{opfqrec}, respectively. The upper limit and lower limit of voltage magnitude at each bus are given in \eqref{opff}, where the voltage limits are presented in their square form.
Equations \eqref{opflimrec} and \eqref{opflimsend} present the linearized formulation for thermal line limits, where $\varepsilon_l$ is an auxiliary parameter, which is dependent on the power factor of load as calculated in \citep{Manshadi2015ResilientMicrogrids}. Here, the absolute values of the real and reactive power flow of line $l$ at time $t$ are modeled. To procure the absolute values, constraints \eqref{opfrec_p}-\eqref{opfsend_q} are employed with a set of auxiliary variables.

\begin{subequations} \label{opf}
\begin{alignat}{3}
&\underset{C_{g,P}^{'u,t}}{\text{max}} \sum_{t \in \mathcal{T}}^{}\sum_{g \in \mathcal{G}_B}^{}(\sum_{i \in \mathcal{I}_g}^{}\lambda_i^{P,t}P_g^t+\lambda_i^{Q,t}Q_g^t-C_g^QQ_g^t-\sum_{u \in \mathcal{U}_g}^{ }C_g^{P,u}P_g^{u,t})\label{opfUL}\\%(Q_g^{+,t}+Q_g^{-,t})
%&\text{s.t.}\nonumber\\
&\text{s.t.}\hspace{1cm}
\underline{\alpha}C_g^{P,u} \leq C_{g,P}^{'u,t} \leq \overline{\alpha}C_g^{P,u}\label{opf_lim_bid}\\
&\underset{P_g^{u,t},Q_g^{u,t}}{\text{min}}\sum_{t \in \mathcal{T}}^{}(\sum_{g\in  \mathcal{G}_B}^{ }\sum_{u \in \mathcal{U}_g}^{ }{C_{g,P}^{'u,t}}P_g^{u,t}+\sum_{g\in  \mathcal{G}_{-B}}^{ }\sum_{u \in \mathcal{U}_g}^{ }{C_g^{P,u}}P_g^{u,t}+\sum_{g\in  \mathcal{G}}^{ }C_g^QQ_g^t-\sum_{d \in \mathcal{D}}\sum_{z \in \mathcal{Z}_d}^{}W_d^{z,t}P_d^{z,t})\label{opfLL}\\%(Q_g^{+,t}+Q_g^{-,t})
%\sum_{s \in \mathcal{S}_g}^{ }(\sum_{g\in  \mathcal{G}_B}^{ }C_g^sP_{g}^{s}+\sum_{s \in \mathcal{G}_{-B}}^{}{C_g^s}'P_{g}^{s})
%&\text{s.t.}\nonumber\\
&\text{s.t.}\hspace{0.5cm}\sum_{u \in \mathcal{U}_g}^{ } P_{g}^{u,t}=P_g^t&\hspace{-3cm}:\lambda_{P_{g,t}} \label{opfpseg}\\
% & Q_g^t=Q_g^{+,t}-Q_g^{-,t} \hspace{3.5cm}:\lambda_{Q_{g,t}} \label{opfQ}\\
% &\sum_{u \in \mathcal{U}_g}^{ } Q_{g}^{u,t}=Q_g^t\hspace{4.4cm}:\lambda_{Q_{g,t}} \label{opfqseg}\\
&-R_g^d \leq P_g^t-P_g^{t-1} \leq R_g^u\hspace{1.2cm}\forall t \in \mathcal{T}\setminus{1}&\hspace{-3cm}:\underline{\mu}_{P_g^{r,t}},\overline{\mu}_{P_g^{r,t}} \label{opf_ramp_t}\\
&-R_g^d \leq P_g^t-P_g^{|\mathcal{T}|} \leq R_g^u\hspace{1.6cm}\forall t=1&\hspace{-3cm}:\underline{\mu}_{P_g^{r,t}},\overline{\mu}_{P_g^{r,t}}\label{opf_ramp_1}\\
&0\leq  P_g^{u,t}\leq \overline{P}_g^u&\hspace{-3cm}:\underline{\mu}_{P_g^{u,t}},\overline{\mu}_{P_g^{u,t}}\label{opflimpseg}\\
&\underline{P}_g\leq  P_g^t \leq \overline{P}_g&\hspace{-3cm}
:\underline{\mu}_{g}^{P,t},\overline{\mu}_{g}^{P,t}\label{opfpg}\\
&\underline{Q}_g+\sum_{u \in \mathcal{U}_g}\rho_g^{+,u}P_g^{u,t} \leq  Q_g^t\leq \overline{Q}_g+\sum_{u \in \mathcal{U}_g}\rho_g^{-,u}P_g^{u,t}&\hspace{-3cm} :\underline{\mu}_{Q^t_g},\overline{\mu}_{Q^t_g}\label{opfqg_p}\\
&0 \leq P_d^{z,t} \leq \chi_d^z p_d^t \hspace{2.9cm}\forall d \in \mathcal{D}\setminus{\mathcal{D}_c}&\hspace{-3cm}:\underline{\mu}_{P_d^{z,t}},\overline{\mu}_{P_d^{z,t}}\label{willingness_p}\\
&\sum_{z \in\mathcal{Z}_d} P_d^{z,t} = p_d^t\hspace{3.5cm} \forall d \in \mathcal{D}\setminus{\mathcal{D}_c} &\hspace{-3cm}:\lambda_{P_d^{fi,t}}\label{opf_fixed}\\
&\left\{\begin{matrix}
\sum_{t \in \mathcal{T}}^{} P_d^{z,t} = \chi_d^z\sum_{t \in \mathcal{T}}^{}p_d^t\hspace{1.3cm} \forall d \in \mathcal{D}_c, z \in \mathcal{Z}_d^{sh}\hspace{4.9cm}:\lambda_{P_d^f}
\\ 
 P_d^{z,t} = \chi_d^zp_d^t\hspace{3.3cm} \forall d \in \mathcal{D}_c, z \in \mathcal{Z}_d\setminus{\mathcal{Z}_d^{sh}}\hspace{4.5cm}:\lambda_{P_d^{f,t}}
\end{matrix}\right.\label{opf_flex}\\
&\sum_{g \in \mathcal{G}_i } P_g^t+\sum_{l\in \mathcal{LT}_i }P_l^{r,t}=\sum_{l\in \mathcal{LF}_i}P_l^{s,t}+\sum_{z \in \mathcal{Z}_d}^{}P_d^{z,t}&\hspace{-3cm}:\lambda_i^{P,t} \label{opfpbus}\\
&\sum_{g \in \mathcal{G}_i } Q_g^t+\sum_{l\in \mathcal{LT}_i }Q_l^{r,t}=\sum_{l\in \mathcal{LF}_i}Q_l^{s,t}+q_d^t  &\hspace{-3cm}:\lambda_i^{Q,t} \label{opfqbus}\\
&P_l^{s,t}=\sum_{i \in \mathcal{F}_l , j \in \mathcal{T}_l}^{ } (g_{l}({e_i^t}^2+{f_i^t}^2)-g_{l}(e_i^te_j^t+f_i^tf_j^t)-b_{l}(e_j^tf_i^t-e_i^tf_j^t))&\hspace{-3cm}:\lambda_l^{s,P,t}\label{opfpsend}\\
&P_l^{r,t}=\sum_{i \in \mathcal{F}_l , j \in \mathcal{T}_l}^{ }(-g_{l}({e_j^t}^2+{f_j^t}^2)+g_{l}(e_i^te_j^t+f_i^tf_j^t)+b_{l}(e_i^tf_j^t-e_j^tf_i^t))&\hspace{-3cm}:\lambda_l^{r,P,t}\label{opfprec}\\
&Q_l^{s,t}=\sum_{i \in \mathcal{F}_l , j \in \mathcal{T}_l}^{ }(-(b_{l}+\frac{b_{l}^{sh}}{2})({e_i^t}^2+{f_i^t}^2)+b_{l}(e_i^te_j^t+f_i^tf_j^t)-g_{l}(e_j^tf_i^t-e_i^tf_j^t))&\hspace{-3cm}:\lambda_l^{s,Q,t}\label{opfqsend}
\\
&Q_l^{r,t}=\sum_{i \in \mathcal{F}_l , j \in \mathcal{T}_l}^{ }((b_{l}+\frac{b_{l}^{sh}}{2})({e_j^t}^2+{f_j^t}^2)-b_{l}(e_i^te_j^t+f_i^tf_j^t)+g_{l}(e_i^tf_j^t-e_j^tf_i^t))&\hspace{-3cm}:\lambda_l^{r,Q,t}\label{opfqrec} \\
&({\underline{V}_i})^2\leq  {e_i^t}^2+{f_i^t}^2\leq ({\overline{V}_i})^2 &\hspace{-3cm}: \underline{\mu}_{v_i^t},\overline{\mu}_{v_i^t}\label{opff}\\
&(P_{l,+}^{r,t}+P_{l,-}^{r,t})+\varepsilon_l (Q_{l,+}^{r,t}+Q_{l,-}^{r,t})\leq \overline{S_l}&\hspace{-3cm}:\overline{\mu}_l^{r,t}\label{opflimrec}\\
&(P_{l,+}^{s,t}+P_{l,-}^{s,t})+\varepsilon_l (Q_{l,+}^{s,t}+Q_{l,-}^{s,t})\leq \overline{S_l}&\hspace{-3cm}:\overline{\mu}_l^{s,t}\label{opflimsend}\\
&P_{l,+}^{r,t}-P_{l,-}^{r,t}=P_l^{r,t}&\hspace{-3cm}:\lambda_{P_l^{r,t}} \label{opfrec_p}\\ &P_{l,+}^{s,t}-P_{l,-}^{s,t}=P_l^{s,t}&\hspace{-3cm}:\lambda_{P_l^{s,t}}\label{opfsend_p}\\
& Q_{l,+}^{r,t}-Q_{l,-}^{r,t}=Q_l^{r,t}&\hspace{-3cm}:\lambda_{Q_l^{r,t}}\label{opfrec_q}\\
&Q_{l,+}^{s,t}-Q_{l,-}^{s,t}=Q_l^{s,t}&\hspace{-3cm}:\lambda_{Q_l^{s,t}}\label{opfsend_q}
\end{alignat}
\end{subequations}
The problem formulation presented in \eqref{opf} is a bi-level optimization problem. The lower-level problem is a non-convex quadratic optimization problem. The non-convexity of the lower-level problem is due to the bi-linear terms in the branch flow equality constraints given in \eqref{opfpsend}-\eqref{opff}. %In this case, because of the non-convexity of the lower-level problem, the duality gap exists. 

%To reformulate the bi-level strategic bidding optimization problem as a single-level optimization problem, the dual form of the lower-level problem is presented. Because of the duality gap in the lower-level problem, the dual variables are not exact. %The bidding strategy of each generation unit depends on LMP at the generator bus, therefore the bidding strategy is not reliable. 
To find a tractable solution to the presented bi-level optimization problem presented in this section, it should be converted to a single-level optimization problem which requires the MCP to be presentable in a closed-form. The current non-convex form cannot be presented in the closed-form. Thus, a solution method to present the dual formulation of the relaxed form of the non-convex problem in the lower-level is shown in the next section. 

%\vspace{-0.1cm}
\section {Solution Methodology}\label{sol}
\subsection{Overview}
A convex relaxation scheme is first presented in this section to facilitate the formation of the closed-form representation for the MCP given in \eqref{opfLL}-\eqref{opflimsend}.
The MCP is reformulated using the SOCP relaxation method. Then, the dual problem associated with the SOCP problem is introduced in \eqref{dual_opf}. The duality gap for the SOCP problem is zero if the conic constraints are strictly feasible \citep{lobo1998applications}. Next, the primal-dual pair constraints are employed to
procure the closed-form of the MCP, as presented in \eqref{closed}.
The procured closed-form enables the reformulation of the original strategic bidding problem in \eqref{opf} into a single-level problem. However, several nonlinear terms will remain in the single-level optimization problem. The final step is to present valid constraints to reformulate the nonlinear terms into equivalent linear terms.

\subsection{The SOCP Relaxation of the Market-Clearing Problem} \label{relax}
A set of lifting variables is introduced in \eqref{cc} to relax the nonlinear terms in the lower-level
problem presented in \eqref{opfLL}-\eqref{opflimsend}.

\begin{equation} \label{cc}
c_{ii}^t:={e_i^t}^2+{f_i^t}^2 ; \hspace{2mm} c_{ij}^t:=e_i^te_j^t+f_i^tf_j^t ; \hspace{2mm} s_{ij}^t:=e_j^tf_i^t-f_j^te_i^t%=\textbf{Im}\{V_iV_j^*
\end{equation}
The lifting terms \eqref{cc} are employed to reformulate the original MCP to the form presented in \eqref{Popf}. This relaxation is referred to as the SOCP relaxation of the ACOPF problem. Here, the objective function is the same as the one in the non-convex lower-level problem. The non-convex equality constraints presented in \eqref{opfpsend}-\eqref{opfqrec} are represented in the form given in \eqref{popfpsend}-\eqref{popfqrec} by leveraging the lifting variables introduced in \eqref{cc}. The voltage limits are presented in \eqref{popfvol} leveraging the lifting terms. Relationships between the lifting terms of each pair of buses are given in \eqref{popf_symtric}. The second-order cone constraint \eqref{popfSOC} presents the second-order cone relaxation of the relationship between the lifting terms presented in \eqref{cc}. The rest of the constraints are the same as those in the original MCP.
%\color{black}{Non-convex original ACOPF formulation goes here} \color{black}
%\vspace{-0.25cm}
\begin{subequations} \label{Popf}
\begin{alignat}{3}
&\underset{P_g^{u,t},Q_g^{u,t}}{\text{min}}\sum_{t \in \mathcal{T}}^{}(\sum_{g\in  \mathcal{G}_B}^{ }\sum_{u \in \mathcal{U}_g}^{ }{C_{g,P}^{'u,t}}P_g^{u,t}+\sum_{g\in  \mathcal{G}_{-B}}^{ }\sum_{u \in \mathcal{U}_g}^{ }{C_g^{P,u}}P_g^{u,t}+\sum_{g\in  \mathcal{G}}^{ }C_g^QQ_g^t-\sum_{d \in \mathcal{D}}\sum_{z \in \mathcal{Z}_d}^{}W_d^{z,t}P_d^{z,t})\label{Popfa}\\
&\text{s.t.}\hspace{0.2cm}P_l^{s,t}=\sum_{i \in \mathcal{F}_l , j \in \mathcal{T}_l}^{ } g_{l}c_{ii}^t-g_{l}c_{ij}^t-b_{l}s_{ij}^t&\hspace{-3cm}:\lambda_l^{s,P,t}\label{popfpsend}\\
&P_l^{r,t}=\sum_{i \in \mathcal{F}_l , j \in \mathcal{T}_l}^{ }-g_{l}c_{jj}^t+g_{l}c_{ji}^t+b_{l}s_{ji}^t&\hspace{-3cm}:\lambda_l^{r,P,t}\label{popfprec}\\
&Q_l^{s,t}=\sum_{i \in \mathcal{F}_l , j \in \mathcal{T}_l}^{ }-(b_{l}+\frac{b_l^{sh}}{2})c_{ii}^t+b_{l}c_{ij}^t-g_{l}s_{ij}^t&\hspace{-3cm} :\lambda_l^{s,Q,t}\label{popfqsend}
\\
&Q_l^{r,t}=\sum_{i \in \mathcal{F}_l , j \in \mathcal{T}_l}^{ }(b_{l}+\frac{b_l^{sh}}{2})c_{jj}^t-b_{l}c_{ji}^t+g_{l}s_{ji}^t&\hspace{-3cm} :\lambda_l^{r,Q,t}\label{popfqrec}\\
&({\underline{V}_i})^2\leq  c_{ii}^t \leq ({\overline{V}_i})^2 \hspace{2.4cm}&\hspace{-3cm}: \underline{\mu}_{v_i^t},\overline{\mu}_{v_i^t} \label{popfvol}\\ 
&c_{ij}^t=c_{ji}^t,\hspace{0.3cm} s_{ij}^t=-s_{ji}^t &\hspace{-3cm} :\lambda_{bp}^{c,t}, \lambda_{bp}^{s,t}\label{popf_symtric}\\
&    \begin{Vmatrix}
2c_{ij}^t
\\2s_{ij}^t
\\ 
c_{ii}^t-c_{jj}^t
\end{Vmatrix}\leq c_{ii}^t+c_{jj}^t&\hspace{-3cm}:\mu_{bp}^{c,t},\mu_{bp}^{s,t},\mu_{bp}^{cc,t},\lambda_{bp}^{cone,t}\label{popfSOC}\\
&\eqref{opfpseg}-\eqref{opfqbus}, \eqref{opfrec_p}, \eqref{opflimsend}\nonumber
\end{alignat}
\end{subequations}
Note that dual variables associated with each constraint in the SOCP relaxation problem \eqref{Popf} are indicated at their related constraint following a colon. In the next subsection, the dual form of the SOCP relaxation problem given in \eqref{Popf} is presented.

\subsection{The Dual Form of the Relaxed Market-Clearing Problem}
The dual formulation of the MCP is presented in \eqref{dual_opf}. Note that the strategic bidding variables $C_{g,P}^{'u,t}$ are the summation of binary variables multiplied by various constants (i.e., $C_{g,P}^{'u,t}=\sum_{c=1}^{N_c}\alpha_cI_g^{c,u,t}C_g^{P,u}$), and it is treated as a constant in procuring the dual form of the problem.
The objective of the dual problem is presented in \eqref{dual_opfa}, where it is procured by summation of the multiplying right-hand side constants of the constraints within the primal problem by associated dual variables. The dual constraints associated with real power dispatch of segment $u$ of generation unit $g$ at time $t$ are given in \eqref{dual_opf_pgs} and \eqref{dual_opf_pgs_B}. The dual constraints corresponding with the real and reactive power dispatch of generation unit $g$ at time $t$ are presented in \eqref{dual_opf_pgt}-\eqref{dual_opf_qg}.
The dual constraints associated with delivered power to the shiftable and non-shiftable portions of curtailable load $d$ at time $t$ are presented in \eqref{dual_opf_pdz_f}.
The dual constraint associated with the portion $z$ of delivered power to fixed load $d$ at time $t$ is presented in \eqref{dual_opf_pdz}. 
The dual constraints corresponding with the real and reactive power sending to line $l$ at time $t$ are represented in \eqref{dual_opf_pls} and \eqref{dual_opf_qls}, respectively. The dual constraints corresponding to the real and reactive power received from line $l$ at time $t$ are given in \eqref{dual_opf_plr} and \eqref{dual_opf_qlr}, respectively. The dual constraints associated with the non-negative variables representing the positive and negative part of real and reactive power flow receiving/sending from/to line $l$ at time $t$ are presented in \eqref{dual_opf_p_rec_+}-\eqref{dual_opf_q_send_-}. The dual constraints associated with the lifting terms $c_{ii}^t$, $c_{ij}^t$, $c_{ji}^t$, $s_{ij}^t$, and $s_{ji}^t$ are given in \eqref{dual_opf_cii}-\eqref{dual_opf_sji}, respectively.
The dual conic constraint associated with second-order cone constraint \eqref{popfSOC} is shown in \eqref{dual_opf_con}. The dual variables associated with the primal inequality constraints and the right-hand side of the dual second-order cone are non-negative, as presented in \eqref{dual_opf_con}.%\vspace{-0.3cm}
\begin{subequations} \label{dual_opf}
\begin{alignat}{3}
&\text{max}\sum_{t \in \mathcal{T}}^{}[\sum_{g\in  \mathcal{G} }^{ } -\underline{\mu}_{P_g^{r,t}}R_g^d-\overline{\mu}_{P_g^{r,t}}R_g^u-\overline{Q}_g\overline{\mu}_{Q_g^t}+\underline{Q}_g\underline{\mu}_{Q_g^t}+\underline{\mu}_{g}^{P,t}\underline{P}_g-\overline{\mu}_{g}^{P,t}\overline{P}_g-\nonumber\\
&\hspace{0.2cm}\sum_{u \in \mathcal{U}_g}^{ }\overline{\mu}_{P_g^{u,t}}\overline{P}_g^u-\sum_{l\in  \mathcal{L} }^{ } (\overline{\mu}_{l}^{r,t}+\overline{\mu}_{l}^{s,t})\overline{S_l}-\sum_{d \in \mathcal{D}\setminus{\mathcal{D}_c}}\sum_{z \in \mathcal{Z}_d}^{}\overline{\mu}_{P_d^{z,t}}\chi_d^zp_d^t+\sum_{d \in \mathcal{D}_c}\sum_{z \in \mathcal{Z}_d^{sh}}^{}\lambda_{P_d^f}\chi_d^zp_d^t\nonumber\\
&\hspace{0.2cm}+\sum_{z \in \mathcal{Z}_d\setminus{\mathcal{Z}_d^{sh}}}^{}\lambda_{P_d^{f,t}}\chi_d^zp_d^t+\sum_{d \in \mathcal{D}\setminus{D}_f}\lambda_{P_d^{fi,t}}p_d^t+\sum_{i\in  \mathcal{N} }^{ } (\sum_{d \in \mathcal{D}_i}^{}\lambda_{i}^{Q,t}q_d^t+\underline{\mu}_{v_i^t}\underline{V}_i^2-\overline{\mu}_{v_i^t}\overline{V}_i^2)]\label{dual_opfa}\\
&\text{s.t.}\nonumber\\
&\underline{\mu}_{P_g^{u,t}}-\overline{\mu}_{P_g^{u,t}}+\lambda_{P_{g,t}}-\rho_g^{+,u}\underline{\mu}_{Q_g^t}+\rho_g^{-,u}\overline{\mu}_{Q_g^t}= C_g^{P,u}\hspace{0.6cm}\forall g \in \mathcal{G}_{-B} \hspace{0.5cm}:P_g^{u,t}\label{dual_opf_pgs}\\
&\underline{\mu}_{P_g^{u,t}}-\overline{\mu}_{P_g^{u,t}}+\lambda_{P_{g,t}}-\rho_g^{+,u}\underline{\mu}_{Q_g^t}+\rho_g^{-,u}\overline{\mu}_{Q_g^t}= C_{g,P}^{'u,t}\hspace{0.6cm}\forall g \in \mathcal{G}_{B}\hspace{0.7cm}:P_g^{u,t}\label{dual_opf_pgs_B}\\
&-\lambda_{P_{g,t}}+\underline{\mu}_{P_g^{r,t}}-\underline{\mu}_{P_g^{r,t+1}}-\overline{\mu}_{P_g^{r,t}}+\overline{\mu}_{P_g^{r,t+1}}+\underline{\mu}_g^{P,t}-\overline{\mu}_g^{P,t}\nonumber\\
&\hspace{0.2cm}+\sum_{i\in \mathcal{I}_g}^{} \lambda_i^{P,t}= 0 \hspace{5.6cm}\forall t \in \mathcal{T}\setminus{\mathcal{|T|}}\hspace{1.0cm} : P_g^t\label{dual_opf_pgt}\\
&-\lambda_{P_{g,t}}+\underline{\mu}_{P_g^{r,t}}-\underline{\mu}_{P_g^{r,1}}-\overline{\mu}_{P_g^{r,t}}+\overline{\mu}_{P_g^{r,1}}+\underline{\mu}_g^{P,t}-\overline{\mu}_g^{P,t}\nonumber\\
&\hspace{0.2cm}+\sum_{i\in \mathcal{I}_g}^{} \lambda_i^{P,t}= 0 \hspace{5.7cm}\forall t =\mathcal{|T|}\hspace{1.7cm} : P_g^t\label{dual_opf_pg24}
\\
&\underline{\mu}_{Q_g^t}-\overline{\mu}_{Q_g^t}+\sum_{i\in \mathcal{I}_g}^{}\lambda_i^{Q,t}=C_g^Q\hspace{6.8cm}: Q_g^t\label{dual_opf_qg}\\
&\left\{\begin{matrix}\underline{\mu}_{P_d^{z,t}}+\lambda_{P_d^f}-\sum_{i \in \mathcal{I}_d}^{}\lambda_i^{P,t}=-W_d^{z,t}\hspace{0.4cm}\forall z \in \mathcal{Z}_d^{sh}
\\ 
\underline{\mu}_{P_d^{z,t}}+\lambda_{P_d^{f,t}}-\sum_{i \in \mathcal{I}_d}^{}\lambda_i^{P,t}=-W_d^{z,t}\hspace{0.4cm}\forall z \in \mathcal{Z}\setminus{\mathcal{Z}_d^{sh}}
\end{matrix}\right.\hspace{0.3cm},\forall d \in \mathcal{D}_c\hspace{0.1cm}:P_d^{z,t}\label{dual_opf_pdz_f}\\
&\underline{\mu}_{P_d^{z,t}}-\overline{\mu}_{P_d^{z,t}}+\lambda_{P_d^{fi,t}}-\sum_{i \in \mathcal{I}_d}^{}\lambda_i^{P,t}=-W_d^{z,t}\hspace{1.1cm}\forall d \in \mathcal{D}\setminus{\mathcal{D}_c}\hspace{1.0cm}:P_d^{z,t}\label{dual_opf_pdz}\\
&\sum_{i \in \mathcal{F}_l }^{ }(-\lambda_i^{P,t})+\lambda_l^{s,P,t}-\lambda_{P_l^{s,t}}=0\hspace{6cm} : P_l^{s,t}\label{dual_opf_pls}
\\ 
&\sum_{i \in \mathcal{T}_l }^{ } \lambda_i^{P,t}+\lambda_l^{r,P,t}-\lambda_{P_l^{r,t}}=0\hspace{6.6cm}: P_l^{r,t}\label{dual_opf_plr}
\\
&\sum _{i \in \mathcal{F}_l }^{  }(-\lambda_i^{Q,t})+\lambda_l^{s,Q,t}-\lambda_{Q_l^{s,t}}=0\hspace{5.9cm} : Q_l^{s,t}\label{dual_opf_qls}
\\ 
&\sum_{i \in \mathcal{T}_l }^{ }\lambda_i^{Q,t}+\lambda_l^{r,Q,t}-\lambda_{Q_l^{r,t}}=0\hspace{6.6cm} : Q_l^{r,t}\label{dual_opf_qlr}
\\
&\lambda_{P_l^{r,t}}-\overline{\mu}_l^{r,t} \leq 0 \hspace{8.5cm}:P_{l,+}^{r,t}\label{dual_opf_p_rec_+}\\
& -\lambda_{P_l^{r,t}}-\overline{\mu}_l^{r,t} \leq 0\hspace{8.0cm}:P_{l,-}^{r,t}\label{dual_opf_p_rec_-}\\
& \lambda_{Q_l^{r,t}}-\epsilon_l\overline{\mu}_l^{r,t} \leq 0\hspace{8.2cm}:Q_{l,+}^{r,t}\label{dual_opf_q_rec_+}\\
& -\lambda_{Q_l^{r,t}}-\epsilon_l\overline{\mu}_l^{r,t} \leq 0 \hspace{7.8cm}:Q_{l,-}^{r,t}\label{dual_opf_q_rec_-}\\
& \lambda_{P_l^{s,t}}-\overline{\mu}_l^{s,t} \leq 0 \hspace{8.5cm}:P_{l,+}^{s,t}\label{dual_opf_p_send_+}\\
& -\lambda_{P_l^{s,t}}-\overline{\mu}_l^{s,t} \leq 0\hspace{8.1cm}:P_{l,-}^{s,t}\label{dual_opf_p_send_-}\\
& \lambda_{Q_l^{s,t}}-\epsilon_l\overline{\mu}_l^{s,t} \leq 0\hspace{8.2cm}:Q_{l,+}^{s,t}\label{dual_opf_q_send_+}\\
& -\lambda_{Q_l^{s,t}}-\epsilon_l\overline{\mu}_l^{s,t} \leq 0 \hspace{8cm}:Q_{l,-}^{s,t}\label{dual_opf_q_send_-}\\
&\sum_{l \in \mathcal{LF}_i}^{ 
}(-g_{l}\lambda_l^{s,P,t}+(b_{l}+\frac{b_{l}^{sh}}{2})\lambda_l^{s,Q,t})+\sum_{l \in \mathcal{LT}_i}^{ } (g_{l}\lambda_l^{r,P,t}-(b_{l}+\frac{b_{l}^{sh}}{2})\lambda_l^{r,Q,t})+\underline{\mu}_{v_i^t}-\overline{\mu}_{v_i^t}\nonumber\\
&\hspace{0.2cm}+\sum_{bp \in \mathcal{BPF}_i}^{ }(\mu_{bp}^{cc,t}+\lambda_{bp}^{cone,t})+\sum_{bp \in \mathcal{BPT}_i}^{}(-\mu_{bp}^{cc,t}+\lambda_{bp}^{cone,t})=0\hspace{1.6cm}: c_{ii}^t\label{dual_opf_cii}
\\
&\sum_{l \in \mathcal{BP}_{l}}^{ }(g_{l}\lambda_l^{s,P,t}-b_{l}\lambda_l^{s,Q,t})+\lambda_{bp,t}^{c}+2\mu_{bp}^{c,t}=0\hspace{0.5cm}\forall i\in \mathcal{F}_{bp},j\in \mathcal{T}_{bp}\hspace{0.8cm} : c_{ij}^t \label{dual_opf_cij}
\\
&\sum_{l \in \mathcal{BP}_{l}}^{ }(-g_{l}\lambda_l^{r,P,t}+b_{l}\lambda_l^{r,Q,t})-\lambda_{bp}^{c,t}=0\hspace{1.7cm}\forall i\in \mathcal{F}_{bp},j\in \mathcal{T}_{bp}\hspace{0.7cm} : c_{ji}^t\label{dual_opf_cji}
\\
&\sum_{l \in \mathcal{BP}_{l}}^{ }(b_{l}\lambda_l^{s,P,t}+g_{l}\lambda_l^{s,Q,t})+\lambda_{bp}^{s,t}+2\mu_{bp}^{s,t}=0\hspace{0.8cm}\forall i\in \mathcal{F}_{bp},j\in \mathcal{T}_{bp}\hspace{0.7cm} : s_{ij}^t\label{dual_opf_sij}\\
&\sum_{l \in \mathcal{L}_{l}}^{ }-(b_{l}\lambda_l^{r,P,t}+g_{l}\lambda_l^{r,Q,t})+\lambda_{bp}^{s,t}=0\hspace{2cm}\forall i\in \mathcal{F}_{bp},j\in \mathcal{T}_{bp}\hspace{0.7cm} : s_{ji}^t\label{dual_opf_sji}\\
&\begin{Vmatrix}
\mu_{bp}^{c,t}
\\
\mu_{bp}^{s,t}
\\
\mu_{bp}^{cc,t}
\end{Vmatrix}\leq \lambda_{bp}^{cone,t}\hspace{1cm},\hspace{1cm}\overline{\mu}_{(.)}^{(.)}, \underline{\mu}_{(.)}^{(.)}, \lambda_{bp}^{cone,t}\geq 0\label{dual_opf_con}
% &\overline{\mu}_{(.)}^{(.)}, \underline{\mu}_{(.)}^{(.)}, \lambda_{bp}^{cone,t}\geq 0\nonumber%\label{dual_opf_mups}
\end{alignat}
\end{subequations}
%\vspace{-0.75cm}
 \subsection{The Closed-Form Representation of the Market-Clearing Problem with ACOPF Formulation}\label{close}
 The problem presented in \eqref{Popf} is a convex optimization problem, and Slater's condition is satisfied if primal or dual problems are strictly feasible. Therefore strong duality holds, i.e., carries zero duality gap \citep{boyd2004convex}. Thus, by setting the primal objective function \eqref{Popfa} equal to the dual objective function \eqref{dual_opfa} as shown in \eqref{primal-dual_obj} and adding all primal and dual constraints, the closed-form representation of the MCP is obtained as presented in \eqref{closed}. 
\begin{subequations}
\begin{alignat}{2}
&\sum_{t \in \mathcal{T}}^{}[\sum_{g\in  \mathcal{G}}^{ } -\underline{\mu}_{P_g^{r,t}}R_g^d-\overline{\mu}_{P_g^{r,t}}R_g^u-\sum_{u \in \mathcal{U}_g}^{ }\overline{\mu}_{P_g^{u,t}}\overline{P}_g^u-\overline{Q}_g\overline{\mu}_{Q_g^t}+\underline{Q}_g\underline{\mu}_{Q_g^t}+\underline{\mu}_{g}^{P,t}\underline{P}_g-\nonumber\\
&\overline{\mu}_{g}^{P,t}\overline{P}_g-\sum_{l\in  \mathcal{L} }^{ } (\overline{\mu}_{l}^{r,t}+\overline{\mu}_{l}^{s,t})\overline{S_l}-\sum_{d \in \mathcal{D}\setminus{\mathcal{D}_c}}\sum_{z \in \mathcal{Z}_d}^{}\overline{\mu}_{P_d^{z,t}}\chi_d^zp_d^t+\sum_{d \in \mathcal{D}_c}\sum_{z \in \mathcal{Z}_d^{sh}}^{}\lambda_{P_d^f}\chi_d^zp_d^t+\nonumber\\
&\sum_{z \in \mathcal{Z}_d\setminus{\mathcal{Z}_d^{sh}}}^{}\lambda_{P_d^{f,t}}\chi_d^zp_d^t+\sum_{d \in \mathcal{D}\setminus{D}_f}\lambda_{P_d^{fi,t}}p_d^t+\sum_{i\in  \mathcal{N} }^{ } (\sum_{d \in \mathcal{D}_i}^{}\lambda_{i}^{Q,t}q_d^t+\underline{\mu}_{v_i^t}\underline{V}_i^2-\overline{\mu}_{v_i^t}\overline{V}_i^2)]=\nonumber\\
&\sum_{t \in \mathcal{T}}^{}(\sum_{g\in  \mathcal{G}_B}^{ }\sum_{u \in \mathcal{U}_g}^{ }{C_{g,P}^{'u,t}}P_g^{u,t}+\sum_{g\in  \mathcal{G}}^{ }C_g^QQ_g^t+\sum_{g\in  \mathcal{G}_{-B}}^{ }\sum_{u \in \mathcal{U}_g}^{ }C_g^{P,u}P_g^{u,t}-\sum_{d \in \mathcal{D}}\sum_{z \in \mathcal{Z}_d}^{}W_d^{z,t}P_d^{z,t})\label{primal-dual_obj}\\
&\eqref{opfpseg}-\eqref{opfqbus}, \eqref{opfrec_p}-\eqref{opfsend_q}, \eqref{popfpsend}-\eqref{popfSOC}, 
\eqref{dual_opf_pgs}-\eqref{dual_opf_con}.\nonumber
\end{alignat}\label{closed}
\end{subequations}

The binary-to-continuous variable multiplication in \eqref{primal-dual_obj} (i.e., $C_{g,P}^{'u,t}P_g^{u,t}$) should be presented in an equivalent linear form as shown in \eqref{phi}-\eqref{bin_cont3}. It should be reminded that $C_{g,P}^{'u,t}=\sum_{c=1}^{N_c}\alpha_cI_g^{c,u,t}C_g^{P,u}$ is a summation of binary variables multiplied by associated constants. Here, $\phi_{P,g}^{c,u,t}$ and $\psi_{P,g}^{c,u,t}$ are non-negative continuous variables introduced as auxiliary variables to facilitate linearization.
\begin{subequations}
\begin{alignat}{2}
&\phi_{P,g}^{c,u,t}=I_g^{c,u,t}P_g^{u,t} \hspace{1cm}, \hspace{1cm}I_g^{c,u,t}\in\{0,1\}\label{phi}\\
&0\leq \phi_{P,g}^{c,u,t}\leq \overline{P}_g^uI_g^{c,u,t}\label{bin_cont1}\\
&\phi_{P,g}^{c,u,t}=P_g^{u,t}-\psi_{P,g}^{c,u,t}\label{bin_cont2}\\
&0\leq \psi_{P,g}^{c,u,t}\leq \overline{P}_g^u(1-I_g^{c,u,t})\label{bin_cont3}
\end{alignat}\label{bin_cont}
\end{subequations}

 The first term on the LHS of \eqref{primal-dual_obj} is replaced in \eqref{closed_bidding_obj_lin} by applying the linearization presented in \eqref{bin_cont}. Thus, the primal-dual closed-form of the MCP is reformulated, as shown in \eqref{closed_bidding_lin}.

\begin{subequations}
\begin{alignat}{3}
& \text{L.H.S of }\eqref{primal-dual_obj} =\color{black}\sum_{t \in \mathcal{T}}^{}\sum_{g\in  \mathcal{G}_B}^{ }\sum_{u \in \mathcal{U}_g}^{ }\sum_{c=1}^{N_c}\alpha_c\phi_g^{c,u,t}C_g^{P,u}+\color{black} \nonumber \\
& \text{\hspace{2.5cm} the~rest~of~R.H.S~of~}\eqref{primal-dual_obj}\label{closed_bidding_obj_lin} \\
&\eqref{opfpseg}-\eqref{opfqbus},\eqref{opfrec_p}-\eqref{opfsend_q},\eqref{popfpsend}-\eqref{popfSOC},\eqref{dual_opf_pgs}-\eqref{dual_opf_con}, \eqref{bin_cont1}-\eqref{bin_cont3}.\nonumber
\end{alignat}\label{closed_bidding_lin}
\end{subequations}

\subsection{Tackling the Nonlinearity in Objective Function of the Upper-Level Problem of the Strategic Bidding Problem}\label{tackle}
The MCP presented in \eqref{opf} is reformulated as a set of linear and SOC constraints over continuous and binary variables in \eqref{closed_bidding_lin}. The procured set of constraints can be replaced with the MCP of the original strategic bidding problem given in \eqref{opf}. However, the obtained optimization problem is still non-linear and hard to solve with off-the-shelf mixed-integer conic solvers due to the bi-linear terms in the objective function of the strategic bidding problem \eqref{opfUL}. To relax these bi-linear terms, the non-linear strategic bidding problem \eqref{opf} is represented by an equivalent Mixed-Integer Second-Order Cone Program (MISOCP) problem, as shown in \eqref{bidding}. The details of the process to procure the equivalent form are extensively discussed in the Appendix. 
\begin{subequations}\label{bidding}
\begin{alignat}{2}
&\underset{C_{g,P}^{'u,t}}{\text{max}} \sum_{t \in \mathcal{T}}^{}\sum_{g\in \mathcal{G}_B}(\sum_{u\in \mathcal{U}_g}^{ }\overline{\mu}_{P_g^{u,t}}\overline{P}_g^u+\underline{\mu}_{P_g^{r,t}}R_g^d+\overline{\mu}_{P_g^{r,t}}R_g^u-\underline{\mu}_g^{P,t}\underline{P}_g+\nonumber\\
&\hspace{0.2cm}\overline{\mu}_g^{P,t}\overline{P}_g-\underline{\mu}_{Q_g^t}\underline{Q}_g+\overline{\mu}_{Q_g^t}\overline{Q}_g+\sum_{c=1}^{N_c}\alpha_c\phi_g^{c,u,t}C_g^{P,u}-\sum_{u \in \mathcal{U}_g}^{}C_g^{P,u}P_g^{u,t})\label{reformulate_obj_UL}\\
&\text{s.t.}\hspace{3.9cm}\eqref{closed_bidding_lin}\nonumber
\end{alignat}
\end{subequations}
%\vspace{-0.8cm}
\section{The Special Case of Demand Response Profit Maximization}
%\vspace{-0.1cm}
The profit maximization problem \eqref{bidding} can be generated for the demand side as well. Load service entities with curtailable loads can leverage the proposed strategic bidding approach to maximize their profit or minimize their cost. The objective function for maximizing the profit from the curtailable loads is presented in \eqref{load_bid}. Here, the MCP remains the same as presented in \eqref{closed_bidding_lin}. However, here, curtailable loads submit strategic bids for their WTP while generation units bid marginally.
 %\vspace{-0.2cm}
\begin{align}\label{load_bid}
&\underset{C_{l,P}^{'z,t}}{\text{max}}\sum_{t \in \mathcal{T}}^{}\sum_{d \in \mathcal{D}_c}^{}\sum_{z \in \mathcal{Z}_d}^{}\sum_{i \in \mathcal{I}_d}^{}W_d^{z,t}P_d^{z,t}-\lambda_i^{P,t}P_d^{z,t}+\lambda_i^{Q,t}q_d^t
\end{align}
 %\vspace{-0.9cm}
 
\section {Case Study}
%%\vspace{-0.1cm}
In this section, the performance of the proposed strategic bidding model is evaluated under different scenarios.
Here, CPLEX 12.10 \citep{cplex200711} is employed as the off-the-shelf solver to solve MISOCP programming problems. The presented results are performed on a PC with a Core i7 CPU 4.70GHz processor and 48 GB memory. The base demand of all test cases is set according to the normalized hourly load of California ISO on August 18, 2020.%\vspace{-0.35cm}

%%\vspace{-0.1cm}
\subsection{An Illustrative Example}%\vspace{-0.15cm}
The sample 3-bus test network includes two generation units connected to buses $1$ and $2$, one load connected to bus $3$, and three transmission lines connecting all buses. Table \ref{table_Gen} shows characteristics of generation units, and Table \ref{table_line} shows characteristics of transmission lines.

\begin{table}[h!]\centering
\centering
\caption{Generation Cost Curve of Generation units}
\vspace{0.1cm}
\label{tab:my-table}
\begin{tabular}{cccccc}
\hline \hline
\begin{tabular}[c]{@{}c@{}}Market\\ Participant\end{tabular} & \begin{tabular}[c]{@{}c@{}}Min.\\ Capacity\\ (MW)\end{tabular} & \begin{tabular}[c]{@{}c@{}}Max.\\ Capacity\\ (MW)\end{tabular} & \begin{tabular}[c]{@{}c@{}}$\alpha$\\ (\$/MWh$^2$)\end{tabular} & \begin{tabular}[c]{@{}c@{}}$\beta$\\ (\$/MWh)\end{tabular} & \begin{tabular}[c]{@{}c@{}}$\gamma$\\ (\$)\end{tabular}  \\ \hline \hline
G1                                                           & 0                                                              & 200                                                              & 2.2      & 7.5                         & 4.07                         \\ \hline
G2                                                           & 0                                                            & 250                                                              & 1.16 & 6.7                         & 6.02                         \\ \hline
\end{tabular}\label{table_Gen}
\end{table}
\begin{table}[t!]\centering
\caption{Transmission Network Data}
\label{tab:my-table}
\vspace{0.25cm}
\begin{tabular}{ccccccc}
\hline
Line & \begin{tabular}[c]{@{}c@{}}From\\ Bus\end{tabular} & \begin{tabular}[c]{@{}c@{}}To\\ Bus\end{tabular} & \begin{tabular}[c]{@{}c@{}}Resistance\\ (p.u.)\end{tabular} & \begin{tabular}[c]{@{}c@{}}Reactance\\ (p.u.)\end{tabular} & \begin{tabular}[c]{@{}c@{}}\begin{math}b_l^{sh}\end{math} \\ (p.u.)\end{tabular}& \begin{tabular}[c]{@{}c@{}}Capacity\\ Limit\\ (MW)\end{tabular} \\ \hline\hline
L1   & 2                                                  & 1                                                & $0.00281$                                                          & 0.0281   & 0.00712                                                        & 120                                                               \\ \hline
L2   & 3                                                  & 2                                                & $ 0.00108   $                                                        & $ 0.0108$  & $ 0.01852  $                                                          & 190                                                               \\ \hline
L3   & 1                                                  & 3                                                & $ 0.002.97 $                                                          & $ 0.0297$  & $ 0.00674 $                                                          & 170  \\          \hline                                                
\end{tabular}\label{table_line}
\end{table}

\begin{figure}[h!]
\centering%width=\columnwidth
{\includegraphics[width=0.6\linewidth]{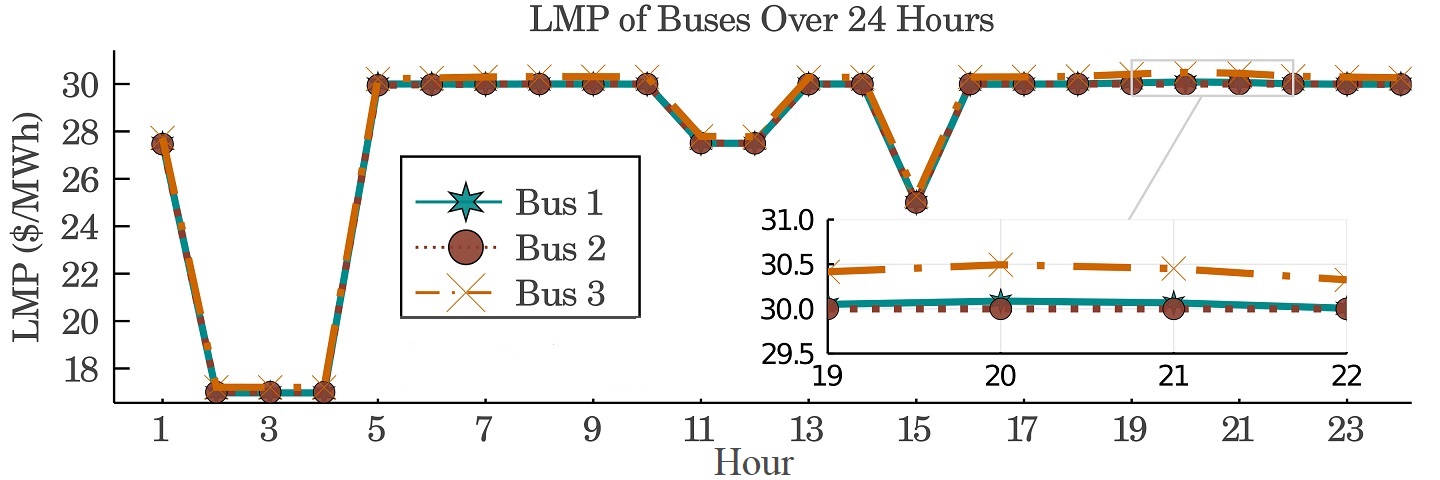}}
\caption{The difference in LMPs of each bus procured by the proposed method}
\label{fig:base_pg}
\end{figure}

Here, the 3-bus system is evaluated under different operational scenarios. The first scenario is the base case of the proposed strategic bidding model; while the ramping constraints of generation units are not enforced, the WTP of the load is much larger than the cost of generation, i.e., the demand is not curtailable. Also, there is no congestion in transmission lines. The proposed strategic bidding model captures the difference between the LMP of different buses in the transmission network because of the loss of transmission lines. Fig. \ref{fig:base_pg} shows the LMP of buses over $24$ hours. An interesting observation here is an increase in the difference between the LMP of the load bus and other buses at hours $19-22$. This is due to the increase in demand, which also leads to an increase in the loss of lines in those hours.

\subsubsection{Comparing strategic bidding models with relaxed ACOPF-based and DCOPF-based market-clearing problem}
 The procured bidding strategies are plugged into the ACOPF-based MCP to illustrate the merit and feasibility of the proposed strategic bidding solution compared with the solution of the strategic bidding problem based on DC-MCP. In the following discussions, the solution procured by plugging the proposed strategic bidding solution into the MCP with the original ACOFP formulation is denoted as RlxAC-AC solution, while the plugged one based on DC-MCP is denoted as DC-AC solution. The clearing prices procured by the proposed strategic bidding problem are in very close proximity to RlxAC-AC prices. This illustrates the effectiveness of the proposed strategic bidding solution as revealed by comparing the third and fifth columns of Table \ref{3bus_dc_SOCP}. Conversely, the strategic bidding based on DC-MCP is different from the DC-AC solution, as revealed by comparing the last row of Table \ref{3bus_dc_SOCP}. Besides, the proposed bidding strategy is more profitable than the one procured by the strategic bidding based on DC-MCP with a $1.25\%$ larger profit margin. The total daily profit of $\$ 26,286.2$ for generation unit $1$ procured by the proposed strategic bidding method is more than that $\$ 25,960$ total daily profit procured by the strategic bidding based on DC-MCP. 
%\vspace{-0.2cm}
\begin{table}[h!] \centering
\caption{Comparing the Solution of the DC and Proposed Strategic Bidding Models for the 3-Bus System at $9$ p.m}
\vspace{0.1cm}
\label{3bus_dc_SOCP}
\begin{tabular}{ccccc} \hline
\multicolumn{1}{c}{\textbf{Method}} & \multicolumn{1}{c}{\textbf{DC}} & \multicolumn{1}{c}{\textbf{RlxAC}} & \multicolumn{1}{c}{\textbf{DC-AC}} & \multicolumn{1}{c}{\textbf{RlxAC-AC}} \\ \hline \hline
$ p_{1}$  [MW]                          & 105  & 66.7 & 105 & 66.7\\ \hline
$ p_{2}$  [MW]                         &   166.7                    &     190  & 166.7  &189.9              \\ \hline
$ \lambda_1^p $                       &  27                     & 30       & 27       & 30        \\ \hline
$ \lambda_2^p $                       &     26.99                  &    30        &27.08   & 30       \\ \hline
$ \lambda_3^p$                      &    \textbf{27.01}                   &   \textbf{30.32}        &\textbf{27.39}  & \textbf{30.31}          \\ \hline
% \begin{tabular}[c]{@{}l@{}}%\gamma_{g1}^{s1}\\ \gamma_{1}^{2}\\ \gamma_{1}^{3}
% $ C_1^{'1}$  \\ $ C_1^{'2}$  \\$ C_1^{'3}$  \end{tabular} &  \begin{tabular}[c]{@{}l@{}}   0\\27\\77      \end{tabular}              &  \begin{tabular}[c]{@{}l@{}}   10\\45\\63      \end{tabular}                        \\ \hline Profit DC_ACOPF= 43589.03, profit DC= 43625.71 , Profit AC
\end{tabular}
\end{table}

\subsubsection{Investigating the impact of transmission lines congestion}
\begin{figure}[h!]
\centering
{\includegraphics[width=0.6\linewidth]{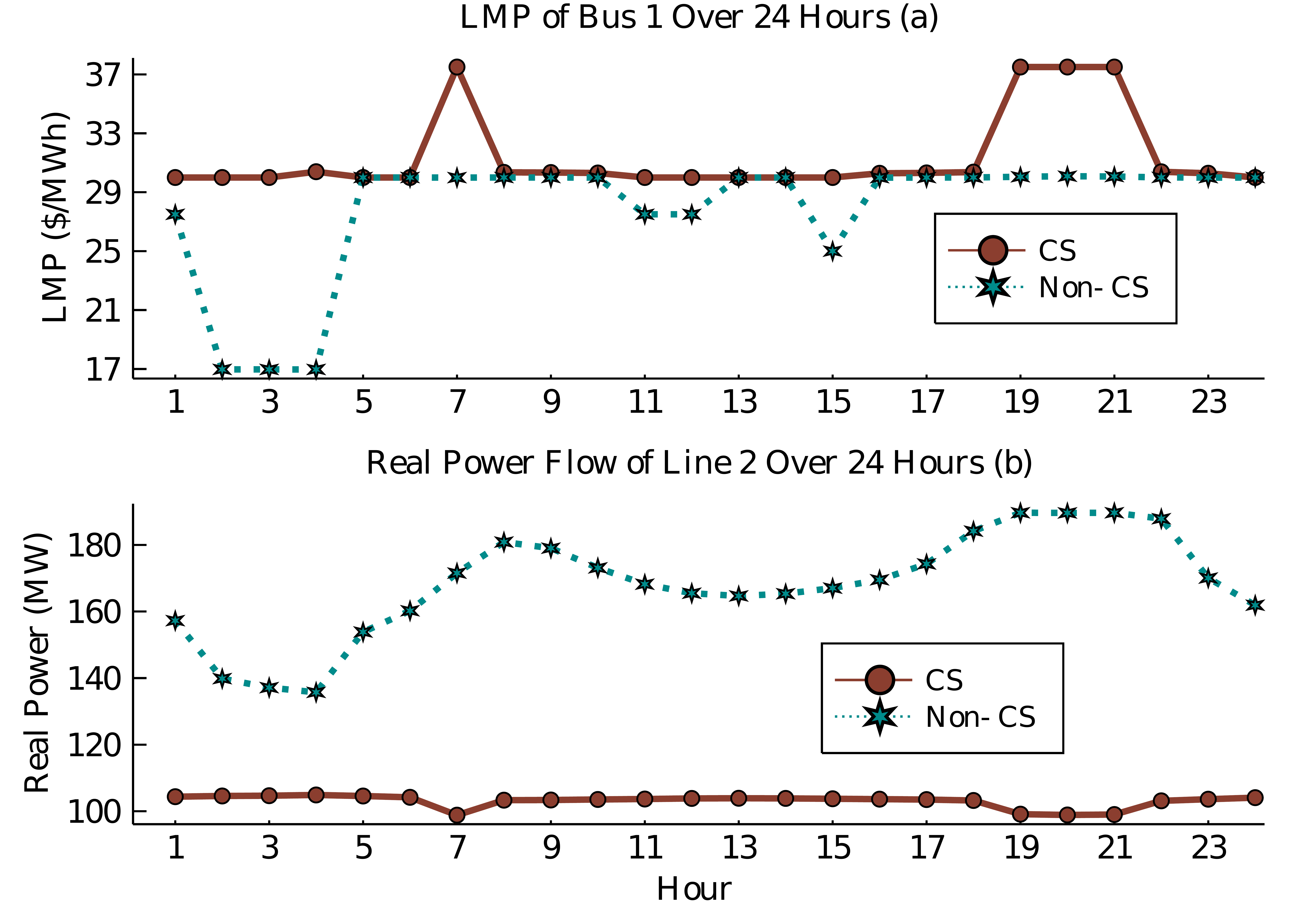}}
%\vspace{-0.8cm}
\caption{The impact of congestion on the real power flow and LMPs}
%\vspace{-0.45cm}
\label{fig:base_cong}
\end{figure}

%Here, the proposed strategic bidding model is evaluated when transmission lines are congested. 
Here, the solution procured by the Non-Congested Scenario (Non-CS) is compared with the one procured by the Congested Scenario (CS). It is assumed that the congestion occurred in line $2$ with a decrease in the thermal capacity. As a result of this congestion, generation unit $1$ responds to net demand increase during peak hours. Thus, the LMP of buses increases during peak hours of $19-21$. Figs. \ref{fig:base_cong}(a) and \ref{fig:base_cong}(b) present the LMP of the load bus and the real power sending from bus $2$ to bus $3$ through line $2$, respectively. Due to the congestion in line $2$ and the increase in the LMP of bus $1$, the profit of generation unit $1$ over $24$ hours in the congested scenario is increased to $\$ 27,874.8$ while the one procured by the non-congested model is $\$ 26,286.2$.

\begin{figure}[h!]
\centering
{\includegraphics[width=0.6\linewidth]{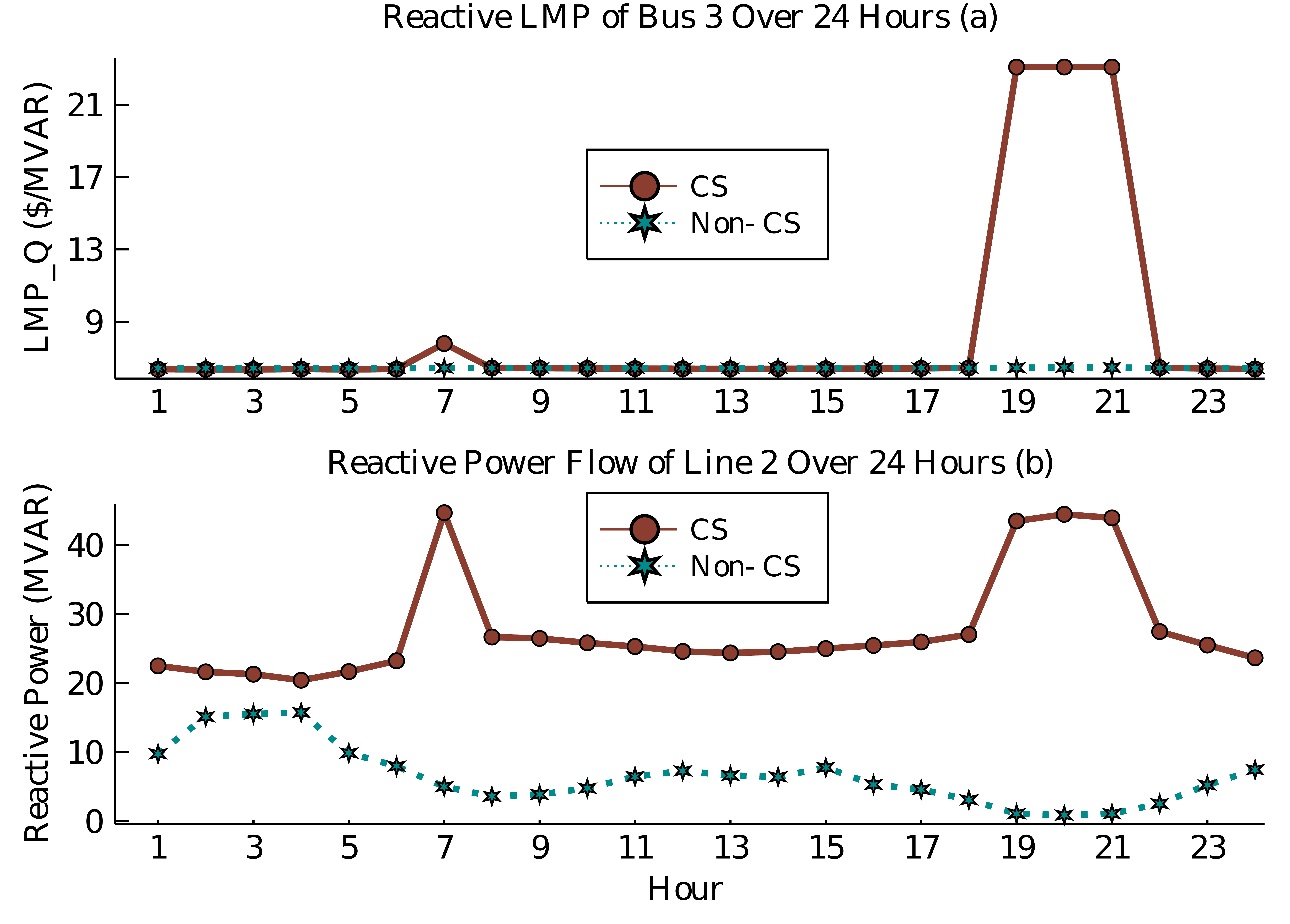}}
%\vspace{-0.4cm}
\caption{The impact of congestion on the reactive power flow and reactive-LMPs}
%\vspace{-0.4cm}
\label{fig:base_cong_q}
\end{figure}

\begin{figure}[h!]
\centering
{\includegraphics[width=0.6\linewidth]{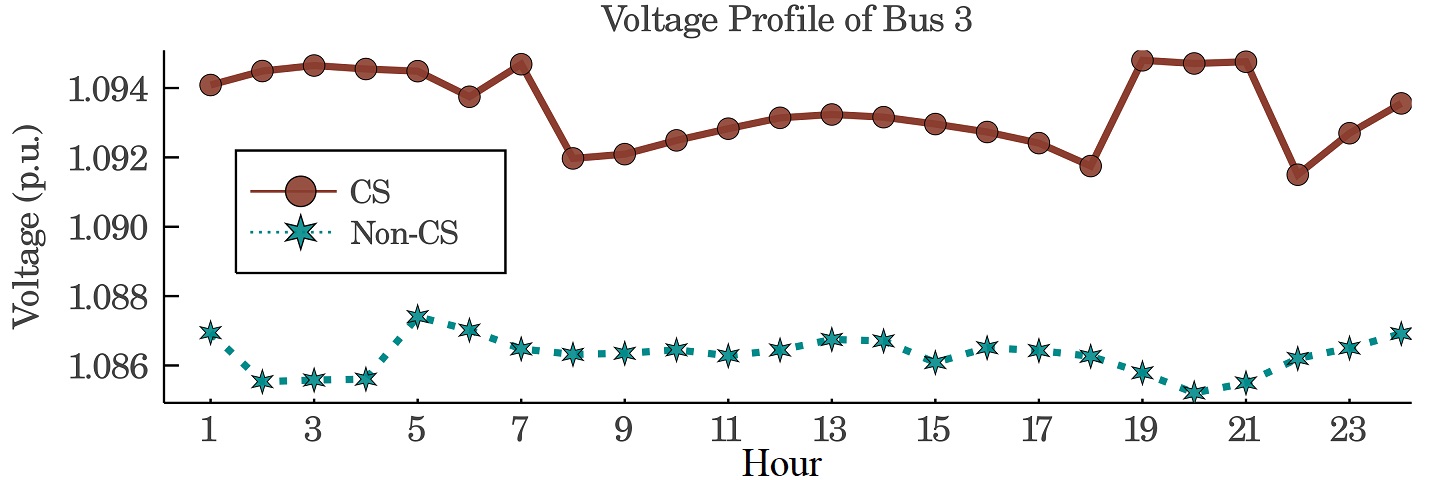}}
%\vspace{-0.5cm}
\caption{The impact of congestion on the voltage profile of bus $3$}
%\vspace{-0.5cm}
\label{fig:base_cong_v_3bus}
\end{figure}

It is interesting to observe the difference in reactive power LMP due to the congestion during peak hours, as shown in Fig. \ref{fig:base_cong_q}(a). Since the reactive power flow passing through line $2$ increases during peak hours to maintain the voltage magnitude of the load bus, which is presented in  Fig. \ref{fig:base_cong_v_3bus} within its acceptable operating range as presented in Fig. \ref{fig:base_cong_q}(b), the real power flow passing through line $2$ decreases during peak hours to hold the line limit constraint.
%\vspace{-0.2cm}
%\vspace{-0.5cm}
%\vspace{-0.4cm}
\vspace{1cm}
 \subsubsection{Investigating the impact of demand response of shiftable loads}
\begin{figure}[h!]
\centering
{\includegraphics[width=0.6\linewidth]{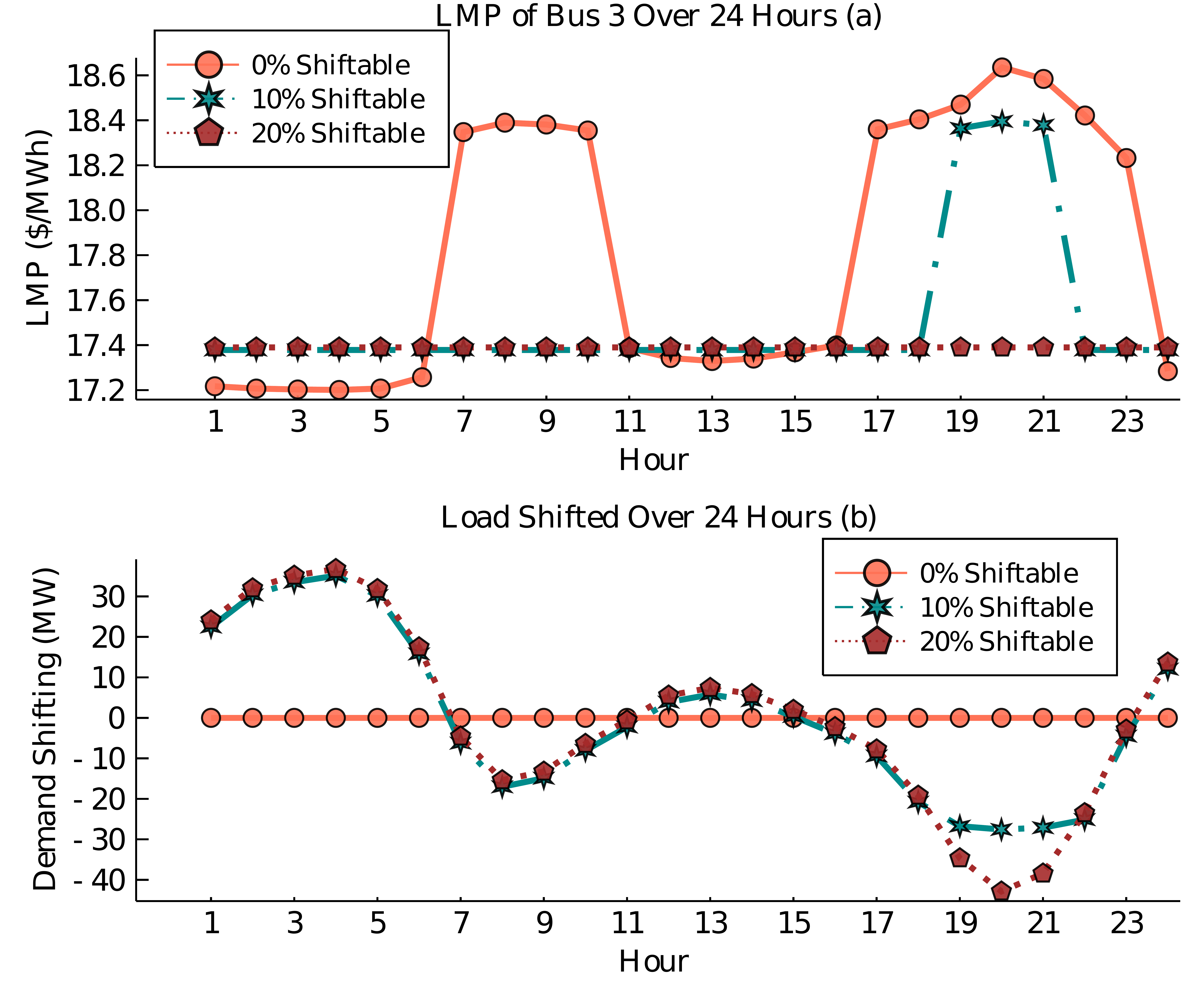}}
%\vspace{-1cm}
\caption{The impact of demand response on the LMPs}
%\vspace{-0.15cm}
\label{fig:shifting_LMP_3}
\end{figure}

 Here, the objective function presented in \eqref{load_bid} is employed to maximize the profit of shiftable loads. Leveraging shiftable loads enables the ISO to serve them during off-peak hours. Thus, increasing the shiftable percentage of the load will mitigate the variation in the LMP of the load bus during the day. Fig. \ref{fig:shifting_LMP_3}(a) presents the LMP of the load bus over $24$ hours when different percentages of the load are shiftable. The LMP of buses when $10\%$ of the load is shiftable is more than when $20\%$ of the load is shiftable for hours $19-21$. Another notable point is that with the $20\%$ shift possibility in the load as shown in  Fig. \ref{fig:shifting_LMP_3}(b), there is no variation in the LMP of the load bus, as shown in Fig. \ref{fig:shifting_LMP_3}(a).
 
 It is interesting to observe that the increase in the demand during off-peak hours as a result of load shifting will lead to a slight increase in the LMP of the load bus during off-peak hours. As shown in Fig. \ref{fig:shifting_LMP_3}(b), the demand is shifted from hours $7-11,17-23$ to hours $1-6,12-15$ which leads to a decrease in the LMP of the load bus at hours $7-11,17-23$ as shown in Fig. \ref{fig:shifting_LMP_3}(a). It is interesting to mention that increasing the shiftable percentage of the load also mitigates the variation in the voltage profile of the load bus.

 %\vspace{-0.2cm}

%  \begin{figure}[h!]
% \centering
% {\includegraphics[width=9.0cm]{figs/congestion_base_volt1.pdf}}
% \caption{The impact of demand response on the voltage profile}
% %\vspace{-0.15cm}
% \label{fig:shifting_volt_1}
% \end{figure}

 \subsubsection{Investigating the impact of ramping limits of generation units}
\begin{figure}[h!]
\centering
{\includegraphics[width=0.6\linewidth]{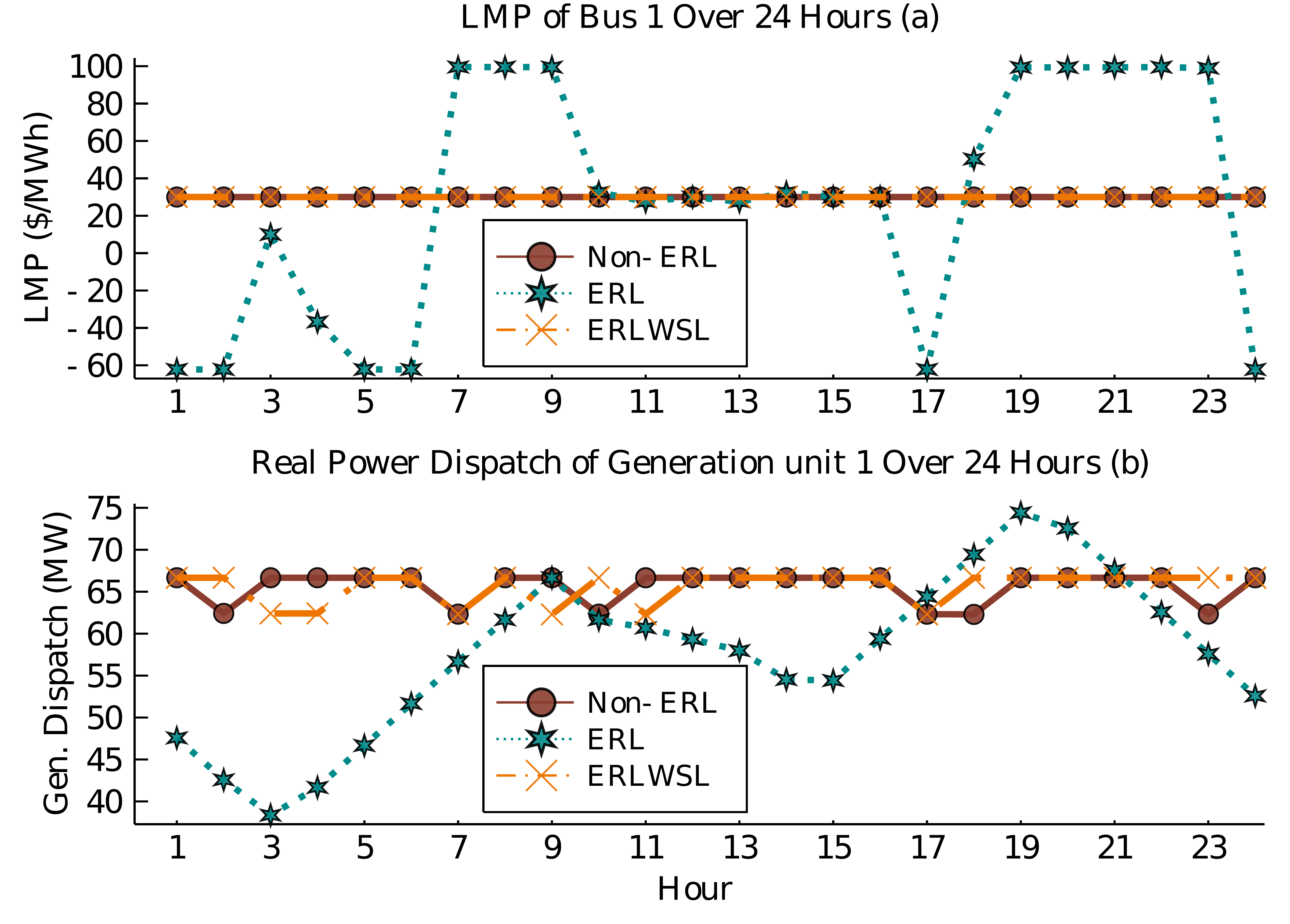}}
%\vspace{-0.75cm}
\caption{The impact of enforcing ramping limits and load shifting}
\label{fig:ramp_base}
\end{figure}

 Here, the impact of ramping constraints of generation units and shiftable loads on the LMP of various buses, real power dispatch of generation units, and the strategic bidding of generation units is discussed. Enforcing ramping limit of generation units may result in load curtailment when the network load is not shiftable; while shifting the load will avoid such an event. 
 Fig. \ref{fig:ramp_base} represents a comparison between the solution procured by the proposed strategic bidding model under three scenarios: with Enforced Ramping Limits (ERL), with Enforced Ramping Limits and Shiftable Loads (ERLWSL), and with Non-Enforced Ramping Limits (Non-ERL). It is interesting to observe that the ramping limit will result in a negative LMP in an hour with low demand, which will increase during the subsequent hours. In other words, the LMP provides an incentive for shiftable demand. Besides, as units will reach their maximum capacity during the peak hours, the LMP of bus $1$ during the peak hours reaches the WTP of curtailable loads as presented in Fig. \ref{fig:ramp_base}(a). However, load shifting mitigated the impact of enforcing ramp limits on the LMP, as shown in Fig. \ref{fig:ramp_base}(b). 

\subsection{IEEE 14-Bus System}
 Here, the modified IEEE 14-bus system is employed when the GENCO, the owner of generation units $1$ and $3$, submits strategic bids. The system contains $5$ generation units and $20$ transmission lines. The proposed strategic bidding solution is compared with the one based on DC-MCP under two scenarios. There is no limitation in the reactive power dispatch of the generation units in the first scenario and the awarded dispatch of the MCP settled the reactive power demand. There is a Limited Reactive Power Support (LRPS) for each generation unit in the second scenario.
 
\begin{figure}[h!]
\centering
{\includegraphics[width=0.6\linewidth]{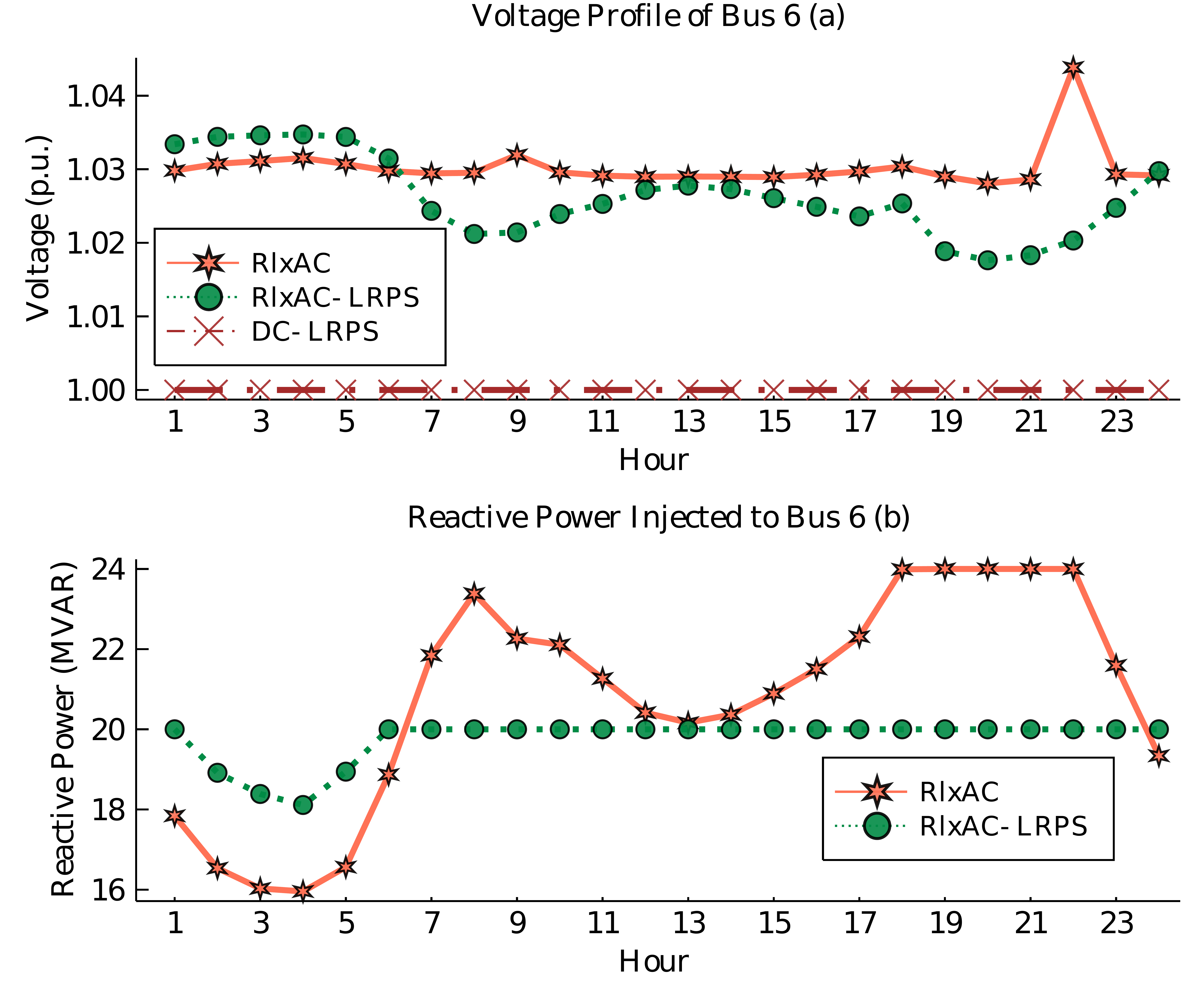}}
%\vspace{-0.9cm}
\caption{The impact of the limitation of reactive power support}
%\vspace{-0.15cm}
\label{fig:LRPS}
\end{figure}

 The strategic bidding model based on DC-MCP cannot correctly model the limitation of reactive power support to modify its bidding strategy accordingly. The limitation of reactive power support leads to a decrease in the voltage magnitude of buses during peak hours. The proposed strategic bidding problem can model the drop in the voltage profile of bus $6$, as shown in Fig. \ref{fig:LRPS} (a). On the contrary, the strategic bidding model based on DC-MCP is unable to model such limitations. As shown in Fig. \ref{fig:LRPS} (b), there is a drop in the reactive power support from the generation unit connected to bus $6$ in case of imposing generation limits. 

Modeling the limitation of reactive power support and considering the reactive power flows result in a more efficient bidding strategy and an increase in the profit of the bidding GENCO. The limited reactive power support will cause an increase in the reactive-LMPs as well as the generation cost. In Table \ref{14bus_dc_SOCP}, the generation cost and total profit of the bidding GENCO procured by the proposed strategic bidding model and the strategic bidding model based on DC-MCP are presented. 
 Besides, the increase in the total profit procured by the proposed strategic bidding model is $52.3\%$ more than the increase in the total profit procured by the strategic bidding model using DC-MCP.  

\begin{table}[h!] \centering
\caption{Comparing the RlxAC-AC and DC-AC Solutions for the 14-Bus System over $24$ Hours}
\label{14bus_dc_SOCP}
\begin{tabular}{ccccc} \hline
\multicolumn{1}{c}{\textbf{Method}} & \multicolumn{1}{c}{\textbf{DC}} & \multicolumn{1}{c}{\textbf{RlxAC}} & \multicolumn{1}{c}{\textbf{DC-LRPS}} & \multicolumn{1}{c}{\textbf{RlxAC-LRPS}} \\ \hline \hline
Total Profit [\$]                          & 17,151.2  & 24,153 & 17,452.3 & 24,611.6\\ \hline
Generation Cost  [\$]                &   65,259.6    &  67,955.7  & 68,249.6   & 68,249.6  \\ \hline
\end{tabular}
\end{table}

%\vspace{-0.4cm}
\subsection{IEEE 118-Bus System}
 In this case, the modified IEEE 118-bus system is considered a test case, and the proposed strategic bidding model is employed to maximize the profit of the GENCO owning generation units $1$ and $25$. The IEEE 118-bus system consists of $54$ generation units, $99$ loads, and $186$ transmission lines.
 The total profit of the GENCO who bids strategically and LMPs procured by the proposed strategic bidding model is in very close proximity to the total profit of the GENCO and LMPs procured by the RlxAC-AC process, and it is more efficient than the one with DC-MCP as manifested in Table \ref{118bus_dc_SOCP}. %The bidding strategy of the proposed method is more efficient than the bidding strategy procured by the strategic bidding model based on DC-MCP as shown in Table \ref{118bus_dc_SOCP}.
 %It is interesting to observe that the proposed strategic bidding model overestimates the profit of the bidding GENCO. However, t
 The strategic bidding model based on DC-MCP underestimates the profit of the bidding GENCO. In comparison, the proposed method has a much smaller gap. Besides, the profit for the RlxAC-AC model is higher than that of the DC-AC one. This illustrates the effectiveness of the proposed approach to procure a meritorious solution. 
 \begin{table}[h!] \centering
 %\vspace{-0.3cm}
\caption{Comparing the Solutions of the Strategic Bidding Models for the 118-Bus System at $7$ p.m}
\label{118bus_dc_SOCP}
\begin{tabular}{cccccc} \hline
\multicolumn{1}{c}{\textbf{Method}} & \multicolumn{1}{c}{\textbf{DC}} & \multicolumn{1}{c}{\textbf{RlxAC}} & \multicolumn{1}{c}{\textbf{DC-AC}} & \multicolumn{1}{c}{\textbf{RlxAC-AC}}&
\multicolumn{1}{c}{\textbf{RlxAC-IB}}\\ \hline \hline
$ p_{1}$  [MW]                          & 66.7  & 33.3 & 66.7 & 33.3 & 33.3\\ \hline
$ p_{25}$  [MW]                         &   85   &  85  & 85  &85  & 85            \\ \hline
$ \lambda_1^p $                       &  20.4 & 21.7       & 19.2       & 21.5 & 21.6        \\ \hline
$ \lambda_{59}^p $                       &     20.6    &    22.2    & 22.1   & 22    & 22   \\ \hline
Total Profit                      &    \textbf{1,047}                   &   \textbf{1,255}        & \textbf{1,168}  & \textbf{1,231}  & \textbf{1,235}        \\ \hline
\end{tabular}
\end{table}
 
 Another interesting observation is the impact of voltage limits on the profit of the GENCO based on AC-MCP. By relaxing the voltage limits from [0.94-1.06] p.u. to [0.9-1.1] p.u., the profit of GENCO is decreased $1.6\%$ as presented in Table \ref{118bus_dc_SOCP}, where the relaxed voltage limits are denoted as RlxAC-IB. %Thus, solving the MCP with AC constraints will has a notable impact on the strategic bidding and the profit of GENCOs. . %When the range of voltage limits increases the reactive power LMP of buses increases to maintain the voltage magnitude of buses within its acceptable operating range. %real and reactive LMPs of buses decreases 

% $ q_{1}$  [MVAR]                         & - & 15  & 76 & 76.7 & 90\\ \hline
% $ q_{25}$  [MVAR]                         & - &  -60   & 16  & 16.2 &     69           \\ \hline
% $ \lambda_1^q $                       &  - & 7.5    & 6.3  & 6.3 &   6.4           \\ \hline
% $ \lambda_{59}^q $                       &     -    &    7.9    & 8  & 8 &  8  \\ \hline
%\vspace{-0.4cm}
\section{Conclusions}
The proposed framework solved the strategic bidding problem with ACOPF-based market-clearing problem. Employing the full ACOPF problem for the market-clearing process allows exploring more opportunities within the strategic bidding problem than the DCOPF-based market-clearing problem.
This problem is a non-convex problem that is hard to solve. This paper proposed an approach to procure a tractable solution that is solvable via off-the-shelf mixed-integer conic solvers (e.g., CPLEX, Gurobi, etc.). To this end, the closed-form primal-dual representation of the relaxed ACOPF-based market-clearing problem is presented to reformulate the bi-level strategic bidding problem as an equivalent single-level mixed-integer cone program. The KKT conditions of the relaxed ACOPF problem and the dual form of the relaxed ACOPF-based market-clearing problem are leveraged to tackle the non-linearity in the objective of the reformulated strategic bidding problem.The presented ACOPF-based strategic bidding enables demand-side management capability with reactive power and ramping constraints of generation units when coupling real and reactive power of generation units. Also, the limitations of reactive power generation dispatch and curtailable loads are explored. 
The performance of the proposed strategic bidding framework is evaluated in the case studies under different scenarios. In the illustrative example, the performance of the proposed model is compared with the conventional strategic bidding with DCOPF-based market-clearing problem. Besides, the merit of the proposed strategic bidding problem, the impact of changing the limits of reactive power support and voltage of buses on the solution of the proposed strategic bidding model, and the profit of the bidding market participant are illustrated. Moreover, the performance of the proposed strategic bidding framework is investigated leveraging the IEEE 14-bus and IEEE 118-bus systems.
It is demonstrated in the case studies that:
\begin{itemize}
   
    \item The clearing prices procured by the proposed strategic bidding problem are in very close proximity to the prices procured by plugging the proposed strategic bidding solution into the market-clearing probles with the original ACOFP formulation.
    
    \item The proposed strategic bidding method procures the solution with more profit since it can model the limitations of reactive power support to update its bids accordingly. In contrast, the strategic bidding problem with DCOPF-based market-clearing problem fails to model such limitations.

    \item The voltage constraint impacts the profit of market participants. It is demonstrated in the IEEE 118-bus system case study that tighter voltage limits will increase the profit of market participants and vice versa. 

\end{itemize}
The discussions of this paper are limited to leveraging the presented solution method for the strategic bidding problem with AC optimal power flow-based market-clearing problem. However, one can apply the proposed method to rendering a tractable solution for bi-level problems that involve a full ACOPF problem formulation in their lower-level problem. Thus, the proposed research paves the way to solve a class of problems on the electricity market and power system resilience.
\section*{Appendix}
Here, our method to tackle the non-linearity of the objective function of the strategic bidding problem is illustrated. Since the strong duality holds for the optimization problem presented in \eqref{Popf}, the KKT conditions are satisfied. Some of the complementary slackness constraints associated with the market-clearing optimization problem given in \eqref{Popf} are presented in \eqref{Slater's}.
%\vspace{-0.1cm}
\begin{subequations}
\begin{alignat}{3}
&\overline{\mu}_{P_g^{u,t}}(P_g^{u,t}-\overline{P}_g^u)=0\Rightarrow \overline{\mu}_{P_g^{u,t}}P_g^{u,t}=\overline{\mu}_{P_g^{u,t}}\overline{P}_g^u  \label{Slater's_mup_over}\\
&\underline{\mu}_{P_g^{u,t}}P_g^{u,t}=0  \label{Slater's_mup_under}\\
&\overline{\mu}_g^{P,t}(P_g^t-\overline{P}_g)=0\Rightarrow \overline{\mu}_g^{P,t}P_g^t=\overline{\mu}_g^{P,t}\overline{P}_g \label{Slater's_mupg_over}\\
&\underline{\mu}_g^{P,t}(\underline{P}_g-P_g^t)=0\Rightarrow \underline{\mu}_g^{P,t}P_g^t=\underline{\mu}_g^{P,t}\underline{P}_g\label{Slater's_mupg_under}\\
&\left\{\begin{matrix}
\underline{\mu}_{P_g^{r,t}}(-R_g^d-P_g^t+P_g^{t-1})=0\hspace{0.2cm}\forall t \in \mathcal{T}\setminus{1}
\\ 
\underline{\mu}_{P_g^{r,t}}(-R_g^d-P_g^t+P_g^{|\mathcal{T}|})=0\hspace{0.6cm}\forall t=1
\end{matrix}\right.\label{slater's_ramping_d}\\
&\left\{\begin{matrix}
\overline{\mu}_{P_g^{r,t}}(P_g^t-P_g^{t-1}-R_g^u)=0\hspace{0.2cm}\forall t \in \mathcal{T}\setminus{1}
\\ 
\overline{\mu}_{P_g^{r,t}}(P_g^t-P_g^{|\mathcal{T}|}-R_g^u)=0\hspace{0.6cm}\forall t=1
\end{matrix}\right.\label{slater's_ramping_u}\\
&\left\{\begin{matrix}
\underline{\mu}_{Q_g^t}(Q_g^t-\underline{Q}_g-\sum_{u \in \mathcal{U}_g}^{}\rho_g^{+,u}P_g^{u,t})=0
\\ 
\overline{\mu}_{Q_g^t}(\overline{Q}_g+\sum_{u \in \mathcal{U}_g}^{}\rho_g^{-,u}P_g^{u,t}-Q_g^t)=0
\end{matrix}\right.\hspace{1cm}\forall g \in \mathcal{G}, u \in  \mathcal{U}_g, t \in  \mathcal{T}\label{slater's_Qg}
%\vspace{-0.3cm}
\end{alignat}\label{Slater's}
\end{subequations}

In the next step, equations \eqref{dual_opf_pgs_B}, \eqref{dual_opf_pgt}, and \eqref{dual_opf_pg24} are leveraged to form the equality presented in \eqref{equality1}. Note that the equation \eqref{equality1} is procured for $t \in \mathcal{T}\setminus{\mathcal{|T|}}$. The similar procedure applies for $t=|T|$.\hspace{-0.5cm}
\begin{align}\label{equality1}
&\underline{\mu}_{P_g^{u,t}}-\overline{\mu}_{P_g^{u,t}}+\underline{\mu}_{P_g^{r,t}}-\underline{\mu}_{P_g^{r,t+1}}-\overline{\mu}_{P_g^{r,t}}+\overline{\mu}_{P_g^{r,t+1}}+\underline{\mu}_g^{P,t}-\overline{\mu}_g^{P,t}-\rho_g^{+,u}\underline{\mu}_{Q_g^t}+\rho_g^{-,u}\overline{\mu}_{Q_g^t}+\sum_{i\in \mathcal{I}_g}^{} \lambda_i^{P,t}= C_{g,P}^{'u,t}\nonumber\\
&\hspace{11cm}\forall g \in \mathcal{G}_B, t \in \mathcal{T}\setminus{\mathcal{|T|}}
\end{align}

Next, by multiplying both sides of equations \eqref{equality1} by $P_g^{u,t}$ and summing up over all segments of generation units and time horizon, the equality given in \eqref{equation3} is procured. $\mathcal{T}\setminus{\mathcal{|T|}}$.%\vspace{-0.2cm}

\begin{align} \label{equation3}
&\sum_{t \in \mathcal{T}}^{}\sum_{i\in \mathcal{I}_g}^{ }\lambda_i^{P,t}P_g^t=\sum_{t \in \mathcal{T}}^{}[\sum_{u\in \mathcal{U}_g}^{ }(-\underline{\mu}_{P_g^{u,t}}P_g^{u,t}+\overline{\mu}_{P_g^{u,t}}P_g^{u,t}-\underline{\mu}_{P_g^{r,t}}P_g^{u,t}+\underline{\mu}_{P_g^{r,t+1}}P_g^{u,t}+\overline{\mu}_{P_g^{r,t}}P_g^{u,t}-\overline{\mu}_{P_g^{r,t+1}}P_g^{u,t}+\nonumber\\
&\rho_g^{+,u}P_g^{u,t}\underline{\mu}_{Q_g^t}-\rho_g^{-,u}P_g^{u,t}\overline{\mu}_{Q_g^t}+C_{g,P}^{'u,t}P_g^{u,t})-\underline{\mu}_g^{P,t}P_g^t+\overline{\mu}_g^{P,t}P_g^t]\hspace{2cm} \forall g \in \mathcal{G}_B
\end{align}

Multiplying the equation \eqref{dual_opf_qg} for all generation units in $\mathcal{G}_B $ by $Q_g^t$ and summing both sides up over time horizon the equation
\eqref{equation_qg} is procured.
\begin{align}\label{equation_qg}
\sum_{t \in \mathcal{T}}^{}\sum_{i\in \mathcal{I}_g}^{}\lambda_i^{Q,t}Q_g^t=\sum_{t \in \mathcal{T}}^{}-\underline{\mu}_{Q_g^t}Q_g^t+\overline{\mu}_{Q_g^t}Q_g^t+C_g^QQ_g^t
\end{align}
Adding both sides of equations \eqref{equation3} and \eqref{equation_qg} and summing both sides up generation units in $\mathcal{G}_B$ the equation \eqref{equation5} is procured.
\begin{align}\label{equation5}
&\sum_{t \in \mathcal{T}}^{}\sum_{g \in \mathcal{G}_B}\sum_{i\in \mathcal{I}_g}^{}\lambda_i^{Q,t}Q_g^t+\lambda_i^{P,t}P_g^t=\sum_{t \in \mathcal{T}}^{}\sum_{g \in \mathcal{G}_B}[-\underline{\mu}_{Q_g^t}Q_g^t+\overline{\mu}_{Q_g^t}Q_g^t+C_g^QQ_g^t+\sum_{u\in \mathcal{U}_g}^{ }(-\underline{\mu}_{P_g^{u,t}}P_g^{u,t}+\overline{\mu}_{P_g^{u,t}}P_g^{u,t}+\nonumber\\
&\rho_g^{+,u}P_g^{u,t}\underline{\mu}_{Q_g^t}-\rho_g^{-,u}P_g^{u,t}\overline{\mu}_{Q_g^t}+C_{g,P}^{'u,t}P_g^{u,t})-\underline{\mu}_{P_g^{r,t}}P_g^{t}+\underline{\mu}_{P_g^{r,t+1}}P_g^{t}+\overline{\mu}_{P_g^{r,t}}P_g^{t}-\overline{\mu}_{P_g^{r,t+1}}P_g^t-\underline{\mu}_g^{P,t}P_g^t+\overline{\mu}_g^{P,t}P_g^t]
\end{align}

Finally, by plugging in complementary slackness equations presented in \eqref{Slater's} into equation \eqref{equation5}, equation \eqref{equation4} is procured. Thus, the bilinear term in \eqref{opfUL} is substituted with their equivalent linear terms, as shown in \eqref{reformulate_obj_UL}.
\begin{align} \label{equation4}
&\sum_{t \in \mathcal{T}}^{}\sum_{g \in \mathcal{G}_B}^{}\sum_{i\in \mathcal{I}_g}^{ }\lambda_i^{Q,t}Q_g^t+\lambda_i^{P,t}P_g^t= \sum_{t \in \mathcal{T}}^{}\sum_{g \in \mathcal{G}_B}^{}[-\underline{\mu}_{Q_g^t}\underline{Q}_g+\overline{\mu}_{Q_g^t}\overline{Q}_g+C_g^QQ_g^t+\sum_{u\in \mathcal{U}_g}^{ }(\overline{\mu}_{P_g^{u,t}}\overline{P}_g^u+\nonumber\\
&\sum_{c=1}^{N_c}\alpha_c\phi_g^{c,u,t}C_g^{P,u})+\underline{\mu}_{P_g^{r,t}}R_g^d+\overline{\mu}_{P_g^{r,t}}R_g^u-\underline{\mu}_g^{P,t}\underline{P}_g+ \overline{\mu}_g^{P,t}\overline{P}_g]
\end{align}
%\vspace{-0.8cm}
\section*{Acknowledgements}
The authors declare that they have no known competing financial interests or personal relationships that could have appeared to influence the work reported in this paper.

\bibliographystyle{elsarticle-num-names} 
\bibliography{cas-refs}

\begin{thebibliography}{27}
\expandafter\ifx\csname natexlab\endcsname\relax\def\natexlab#1{#1}\fi
\providecommand{\url}[1]{\texttt{#1}}
\providecommand{\href}[2]{#2}
\providecommand{\path}[1]{#1}
\providecommand{\DOIprefix}{doi:}
\providecommand{\ArXivprefix}{arXiv:}
\providecommand{\URLprefix}{URL: }
\providecommand{\Pubmedprefix}{pmid:}
\providecommand{\doi}[1]{\href{http://dx.doi.org/#1}{\path{#1}}}
\providecommand{\Pubmed}[1]{\href{pmid:#1}{\path{#1}}}
\providecommand{\bibinfo}[2]{#2}
\ifx\xfnm\relax \def\xfnm[#1]{\unskip,\space#1}\fi
%Type = Article
\bibitem[{Gao et~al.(2015)Gao, Sheble, Hedman, and Yu}]{gao2015optimal}
\bibinfo{author}{F.~Gao}, \bibinfo{author}{G.~B. Sheble},
  \bibinfo{author}{K.~W. Hedman}, \bibinfo{author}{C.-N. Yu},
\newblock \bibinfo{title}{Optimal bidding strategy for {GENCO}s based on
  parametric linear programming considering incomplete information},
\newblock \bibinfo{journal}{International Journal of Electrical Power \& Energy
  Systems} \bibinfo{volume}{66} (\bibinfo{year}{2015})
  \bibinfo{pages}{272--279}.
%Type = Article
\bibitem[{Song et~al.(2002)Song, Liu, and Lawarr{\'e}e}]{song2002nash}
\bibinfo{author}{H.~Song}, \bibinfo{author}{C.-C. Liu},
  \bibinfo{author}{J.~Lawarr{\'e}e},
\newblock \bibinfo{title}{Nash equilibrium bidding strategies in a bilateral
  electricity market},
\newblock \bibinfo{journal}{IEEE transactions on Power Systems}
  \bibinfo{volume}{17} (\bibinfo{year}{2002}) \bibinfo{pages}{73--79}.
%Type = Article
\bibitem[{Li and Shahidehpour(2005)}]{li2005strategic}
\bibinfo{author}{T.~Li}, \bibinfo{author}{M.~Shahidehpour},
\newblock \bibinfo{title}{Strategic bidding of transmission-constrained gencos
  with incomplete information},
\newblock \bibinfo{journal}{IEEE Transactions on power Systems}
  \bibinfo{volume}{20} (\bibinfo{year}{2005}) \bibinfo{pages}{437--447}.
%Type = Article
\bibitem[{Ruiz and Conejo(2009)}]{ruiz2009pool}
\bibinfo{author}{C.~Ruiz}, \bibinfo{author}{A.~J. Conejo},
\newblock \bibinfo{title}{Pool strategy of a producer with endogenous formation
  of locational marginal prices},
\newblock \bibinfo{journal}{IEEE Transactions on Power Systems}
  \bibinfo{volume}{24} (\bibinfo{year}{2009}) \bibinfo{pages}{1855--1866}.
%Type = Article
\bibitem[{Pozo and Contreras(2011)}]{pozo2011finding}
\bibinfo{author}{D.~Pozo}, \bibinfo{author}{J.~Contreras},
\newblock \bibinfo{title}{Finding multiple nash equilibria in pool-based
  markets: A stochastic {EPEC} approach},
\newblock \bibinfo{journal}{IEEE Transactions on Power Systems}
  \bibinfo{volume}{26} (\bibinfo{year}{2011}) \bibinfo{pages}{1744--1752}.
%Type = Article
\bibitem[{Baslis and Bakirtzis(2011)}]{baslis2011mid}
\bibinfo{author}{C.~G. Baslis}, \bibinfo{author}{A.~G. Bakirtzis},
\newblock \bibinfo{title}{Mid-term stochastic scheduling of a price-maker hydro
  producer with pumped storage},
\newblock \bibinfo{journal}{IEEE Transactions on Power Systems}
  \bibinfo{volume}{26} (\bibinfo{year}{2011}) \bibinfo{pages}{1856--1865}.
%Type = Article
\bibitem[{Kazempour et~al.(2011)Kazempour, Conejo, and
  Ruiz}]{kazempour2011strategic}
\bibinfo{author}{S.~J. Kazempour}, \bibinfo{author}{A.~J. Conejo},
  \bibinfo{author}{C.~Ruiz},
\newblock \bibinfo{title}{Strategic generation investment using a
  complementarity approach},
\newblock \bibinfo{journal}{IEEE Transactions on Power Systems}
  \bibinfo{volume}{26} (\bibinfo{year}{2011}) \bibinfo{pages}{940--948}.
%Type = Article
\bibitem[{Hobbs et~al.(2000)Hobbs, Metzler, and Pang}]{hobbs2000strategic}
\bibinfo{author}{B.~F. Hobbs}, \bibinfo{author}{C.~B. Metzler},
  \bibinfo{author}{J.-S. Pang},
\newblock \bibinfo{title}{Strategic gaming analysis for electric power systems:
  An {MPEC} approach},
\newblock \bibinfo{journal}{IEEE transactions on power systems}
  \bibinfo{volume}{15} (\bibinfo{year}{2000}) \bibinfo{pages}{638--645}.
%Type = Article
\bibitem[{Eldridge et~al.(2017)Eldridge, O’Neill, and
  Castillo}]{eldridge2017marginal}
\bibinfo{author}{B.~Eldridge}, \bibinfo{author}{R.~P. O’Neill},
  \bibinfo{author}{A.~Castillo},
\newblock \bibinfo{title}{Marginal loss calculations for the {DCOPF}},
\newblock \bibinfo{journal}{Federal Energy Regulatory Commission, Tech. Rep}
  (\bibinfo{year}{2017}).
%Type = Inproceedings
\bibitem[{Fu and Li(2006)}]{fu2006different}
\bibinfo{author}{Y.~Fu}, \bibinfo{author}{Z.~Li},
\newblock \bibinfo{title}{Different models and properties on {LMP}
  calculations},
\newblock in: \bibinfo{booktitle}{2006 IEEE Power Engineering Society General
  Meeting}, \bibinfo{organization}{IEEE}, \bibinfo{year}{2006}, pp.
  \bibinfo{pages}{1--11}.
%Type = Article
\bibitem[{Sarkar and Khaparde(2009)}]{sarkar2009dcopf}
\bibinfo{author}{V.~Sarkar}, \bibinfo{author}{S.~Khaparde},
\newblock \bibinfo{title}{{DCOPF}-based marginal loss pricing with enhanced
  power flow accuracy by using matrix loss distribution},
\newblock \bibinfo{journal}{IEEE Transactions on Power Systems}
  \bibinfo{volume}{24} (\bibinfo{year}{2009}) \bibinfo{pages}{1435--1445}.
%Type = Article
\bibitem[{Dos~Santos and Diniz(2010)}]{dos2010dynamic}
\bibinfo{author}{T.~N. Dos~Santos}, \bibinfo{author}{A.~L. Diniz},
\newblock \bibinfo{title}{A dynamic piecewise linear model for dc transmission
  losses in optimal scheduling problems},
\newblock \bibinfo{journal}{IEEE Transactions on Power systems}
  \bibinfo{volume}{26} (\bibinfo{year}{2010}) \bibinfo{pages}{508--519}.
%Type = Inproceedings
\bibitem[{Akinbode and Hedman(2013)}]{akinbode2013fictitious}
\bibinfo{author}{O.~W. Akinbode}, \bibinfo{author}{K.~W. Hedman},
\newblock \bibinfo{title}{Fictitious losses in the {DCOPF} with a piecewise
  linear approximation of losses},
\newblock in: \bibinfo{booktitle}{2013 IEEE Power \& Energy Society General
  Meeting}, \bibinfo{organization}{IEEE}, \bibinfo{year}{2013}, pp.
  \bibinfo{pages}{1--5}.
%Type = Article
\bibitem[{Vaishya and Sarkar(2019)}]{vaishya2019accurate}
\bibinfo{author}{S.~Vaishya}, \bibinfo{author}{V.~Sarkar},
\newblock \bibinfo{title}{Accurate loss modelling in the {DCOPF} calculation
  for power markets via static piecewise linear loss approximation based upon
  line loading classification},
\newblock \bibinfo{journal}{Electric Power Systems Research}
  \bibinfo{volume}{170} (\bibinfo{year}{2019}) \bibinfo{pages}{150--157}.
%Type = Article
\bibitem[{Jabr(2008)}]{jabr2008optimal}
\bibinfo{author}{R.~A. Jabr},
\newblock \bibinfo{title}{Optimal power flow using an extended conic quadratic
  formulation},
\newblock \bibinfo{journal}{IEEE transactions on power systems}
  \bibinfo{volume}{23} (\bibinfo{year}{2008}) \bibinfo{pages}{1000--1008}.
%Type = Article
\bibitem[{Lavaei and Low(2011)}]{lavaei2011zero}
\bibinfo{author}{J.~Lavaei}, \bibinfo{author}{S.~H. Low},
\newblock \bibinfo{title}{Zero duality gap in optimal power flow problem},
\newblock \bibinfo{journal}{IEEE Transactions on Power Systems}
  \bibinfo{volume}{27} (\bibinfo{year}{2011}) \bibinfo{pages}{92--107}.
%Type = Inproceedings
\bibitem[{Molzahn and Hiskens(2014)}]{molzahn2014moment}
\bibinfo{author}{D.~K. Molzahn}, \bibinfo{author}{I.~A. Hiskens},
\newblock \bibinfo{title}{Moment-based relaxation of the optimal power flow
  problem},
\newblock in: \bibinfo{booktitle}{2014 Power Systems Computation Conference},
  \bibinfo{organization}{IEEE}, \bibinfo{year}{2014}, pp.
  \bibinfo{pages}{1--7}.
%Type = Article
\bibitem[{Coffrin et~al.(2015)Coffrin, Hijazi, and
  Van~Hentenryck}]{coffrin2015qc}
\bibinfo{author}{C.~Coffrin}, \bibinfo{author}{H.~L. Hijazi},
  \bibinfo{author}{P.~Van~Hentenryck},
\newblock \bibinfo{title}{The qc relaxation: A theoretical and computational
  study on optimal power flow},
\newblock \bibinfo{journal}{IEEE Transactions on Power Systems}
  \bibinfo{volume}{31} (\bibinfo{year}{2015}) \bibinfo{pages}{3008--3018}.
%Type = Article
\bibitem[{Kocuk et~al.(2016)Kocuk, Dey, and Sun}]{kocuk2016strong}
\bibinfo{author}{B.~Kocuk}, \bibinfo{author}{S.~S. Dey}, \bibinfo{author}{X.~A.
  Sun},
\newblock \bibinfo{title}{Strong socp relaxations for the optimal power flow
  problem},
\newblock \bibinfo{journal}{Operations Research} \bibinfo{volume}{64}
  (\bibinfo{year}{2016}) \bibinfo{pages}{1177--1196}.
%Type = Incollection
\bibitem[{Bynum et~al.(2018)Bynum, Castillo, Watson, and
  Laird}]{bynum2018strengthened}
\bibinfo{author}{M.~Bynum}, \bibinfo{author}{A.~Castillo},
  \bibinfo{author}{J.-P. Watson}, \bibinfo{author}{C.~D. Laird},
\newblock \bibinfo{title}{Strengthened socp relaxations for acopf with
  mccormick envelopes and bounds tightening},
\newblock in: \bibinfo{booktitle}{Computer Aided Chemical Engineering},
  volume~\bibinfo{volume}{44}, \bibinfo{publisher}{Elsevier},
  \bibinfo{year}{2018}, pp. \bibinfo{pages}{1555--1560}.
%Type = Inproceedings
\bibitem[{Molzahn and Hiskens(2015)}]{molzahn2015mixed}
\bibinfo{author}{D.~K. Molzahn}, \bibinfo{author}{I.~A. Hiskens},
\newblock \bibinfo{title}{Mixed sdp/socp moment relaxations of the optimal
  power flow problem},
\newblock in: \bibinfo{booktitle}{2015 IEEE Eindhoven PowerTech},
  \bibinfo{organization}{IEEE}, \bibinfo{year}{2015}, pp.
  \bibinfo{pages}{1--6}.
%Type = Article
\bibitem[{Soofi et~al.(2020)Soofi, Manshadi, Liu, and Dai}]{soofi2020socp}
\bibinfo{author}{A.~F. Soofi}, \bibinfo{author}{S.~D. Manshadi},
  \bibinfo{author}{G.~Liu}, \bibinfo{author}{R.~Dai},
\newblock \bibinfo{title}{A socp relaxation for cycle constraints in the
  optimal power flow problem},
\newblock \bibinfo{journal}{IEEE Transactions on Smart Grid}
  (\bibinfo{year}{2020}).
%Type = Article
\bibitem[{Nilsson and Mercurio(1994)}]{nilsson1994synchronous}
\bibinfo{author}{N.~Nilsson}, \bibinfo{author}{J.~Mercurio},
\newblock \bibinfo{title}{Synchronous generator capability curve testing and
  evaluation},
\newblock \bibinfo{journal}{IEEE Transactions on Power Delivery}
  \bibinfo{volume}{9} (\bibinfo{year}{1994}) \bibinfo{pages}{414--424}.
%Type = Article
\bibitem[{Manshadi and Khodayar(2015)}]{Manshadi2015ResilientMicrogrids}
\bibinfo{author}{S.~Manshadi}, \bibinfo{author}{M.~Khodayar},
\newblock \bibinfo{title}{{Resilient operation of multiple energy carrier
  microgrids}},
\newblock \bibinfo{journal}{IEEE Transactions on Smart Grid}
  \bibinfo{volume}{6} (\bibinfo{year}{2015}).
  \DOIprefix\doi{10.1109/TSG.2015.2397318}.
%Type = Article
\bibitem[{Lobo et~al.(1998)Lobo, Vandenberghe, Boyd, and
  Lebret}]{lobo1998applications}
\bibinfo{author}{M.~S. Lobo}, \bibinfo{author}{L.~Vandenberghe},
  \bibinfo{author}{S.~Boyd}, \bibinfo{author}{H.~Lebret},
\newblock \bibinfo{title}{Applications of second-order cone programming},
\newblock \bibinfo{journal}{Linear algebra and its applications}
  \bibinfo{volume}{284} (\bibinfo{year}{1998}) \bibinfo{pages}{193--228}.
%Type = Book
\bibitem[{Boyd et~al.(2004)Boyd, Boyd, and Vandenberghe}]{boyd2004convex}
\bibinfo{author}{S.~Boyd}, \bibinfo{author}{S.~P. Boyd},
  \bibinfo{author}{L.~Vandenberghe}, \bibinfo{title}{Convex optimization},
  \bibinfo{publisher}{Cambridge university press}, \bibinfo{year}{2004}.
%Type = Article
\bibitem[{Cplex(2007)}]{cplex200711}
\bibinfo{author}{I.~Cplex},
\newblock \bibinfo{title}{11.0 user’s manual},
\newblock \bibinfo{journal}{ILOG SA, Gentilly, France}  (\bibinfo{year}{2007})
  \bibinfo{pages}{32}.

\end{thebibliography}





% \end{thebibliography}
\end{document}